\documentstyle[10pt]{article}
\input xypic
\xyoption{all}
\LaTeXdiagrams
\xyoption{2cell}
\UseAllTwocells
\newtheorem{defi}{Definition}
\newtheorem{teo}{Theorem}
\newtheorem{lem}{Lemma}
\newtheorem{prop}{Proposition}
\newtheorem{cor}{Corollary}

\def\be{\begin{equation}}
\def\ee{\end{equation}}
\def\lra{\longrightarrow}
\def\binbar{\overline{\bf Bin}}
\def\fin{Fin_{sk}}
\def\finbar{\overline{\bf Fin}_{sk}}
\def\A{{\bf A}}
\def\B{{\bf B}}
\def\C{{\bf C}}
\def\I{{\bf I}}
\def\Cat{{\bf Cat}}
\def\D{{\bf D}}
\def\Set{{Set}}
\def\TT{{\bf T}}
\def\TTo{{\bf T^{op}}}%
\newcommand{\T}[1]{{\bf T}_{#1}}
\newcommand{\TO}[1]{{\bf T}_{#1}^{op}}
\def\Tassoc{{\bf T}_{Assoc}}
\def\Tmon{{\bf T}_{Mon}}
\def\Tsmon{{\bf T}_{sMon}}
\def\Tsym{{\bf T}_{Sym}}
\def\Tbraid{{\bf T}_{Braid}}
\def\Tsbraid{{\bf T}_{sBraid}}
\def\Tcomm{{\bf T}_{Comm}}
\def\Tscomm{{\bf T}_{sComm}}
\def\Tssym{{\bf T}_{sSym}}%
\def\Tbal{{\bf T}_{Bal}}%
\def\Tsbal{{\bf T}_{sBal}}%
\def\th{{\bf Theories}}
\def\Th{\widetilde{\bf Theories}}
\def\TTh{\widetilde{\bf 2Theories}}
\def\TCatcat{\widetilde{\bf 2Cat/ Cat}}
\def\Catset{\widetilde{\bf Cat /} Set} 
\def\catset{{\bf Cat/}Set}
\def\KP{\otimes^K_0}
\def\TA{{\bf 2Alg}}
\def\TAi{{\bf 2Alg^i}}%
\def\Psir{\Psi_r}%
\def\Psirh{\hat \Psi_r}%
\newcommand{\TAC}[1]{ \TA _{G_{#1}}(\T #1,\C)}
\def\CI{\C ^\I}
\def\CIg{\C ^\I_\gamma}
\def\Psir{\Psi_r}
\def\Psirh{\hat {\Psi}_r}
\newenvironment{examp}{ \stepcounter{examnum} {\bf \noindent Example.\arabic{section}.\arabic{examnum}:}}{$\Box$} 
\newenvironment{remar}{ \stepcounter{remarnum} {\bf \noindent Remark.\arabic{section}.\arabic{remarnum}:}}{$\Box$} 
\newcounter{examnum}[section]
\newcounter{remarnum}[section]
\begin{document}
\title{The Syntax of Coherence}
\author{Noson S. Yanofsky}
\date{September 30, 1999}
\maketitle
\begin{abstract}
\noindent This article tackles categorical coherence 
within a two-dimensional generalization of 
Lawvere's functorial semantics. 2-theories, a 
syntactical way of describing categories with 
structure, are presented. From the
perspective here afforded, many coherence 
results become simple statements about the 
quasi-Yoneda lemma and 2-theory-morphisms. Given 
two 2-theories and a 2-theory-morphism between 
them, we explore the induced relationship 
between the corresponding 2-categories of 
algebras. The strength of the induced 
quasi-adjoints is classified by the strength of 
the 2-theory-morphism. These quasi-adjoints 
reflect the extent to which one structure can be 
replaced by another. A two-dimensional analogue 
of the Kronecker product is defined and 
constructed. This operation allows one to 
generate new coherence laws from old ones.
\end{abstract}
\section{Introduction}

There has been much talk lately about higher-dimensional algebra. One-
dimensional algebra is thought to be about sets with structure.  Many branches of 
mathematics (low-dimensional topology, stable homotopy theory, etc), 
physics (quantum groups, quantum gravity, quantum field theory etc) and computer science (linear 
logic, programming semantics, etc) have made the move from sets with structure 
to categories with structure. This is thought of as two-dimensional algebra. 
One imagines that n-categories with structure would be called n-dimensional 
algebra. This paper is an approach to two-dimensional universal algebra.

Ever since Mac Lane's classic paper \cite{MacLane}, coherence questions have 
played a major role when studying categories with additional structure. 
Coherence deals with the relationship between two operations on a category. 
Whereas when dealing with sets, two operations can either be equal or not 
equal, when dealing with categories, many more options exist. Between any two 
operations on a category, there can be no relation, there can be a morphism, 
there can be an isomorphism, or there can be a unique isomorphism.
Much effort has been exerted to characterize when one structure can 
be replaced by another. These theorems have been proved in an {\it ad hoc} fashion. 
We shall show that many of these theorems can be proven in a universal and 
organic manner.

The formalism that we chose to follow is Lawvere's functorial semantics 
\cite{Lawvere}, \cite{Lawverea}. For each algebraic structure, one constructs a theory $T$ 
whose objects are the natural numbers and whose morphisms $f:n \longrightarrow 
m$ correspond to operations. Composition of morphisms correspond to 
composition and substitution of operations.
Models or algebras of a theory are product 
preserving functors,$F$, from the theory to a category $C$ with finite products. 
So $F(1)$ is an object $c \in C$. $F(n) \sim c^n$ and $F(f:n \longrightarrow 1)$ is an 
$n$-ary operation. Natural transformations between these functors are homomorphisms of the 
structures. Algebras and homomorphisms form a category $Alg(T,C)$. Between 
theories there are theory-morphisms $G:T \longrightarrow T'$. Precomposition 
with such a morphism induces $G^*:Alg(T',C) \longrightarrow Alg(T,C)$. The 
central theory of functorial semantics says that $G^*$ has a left adjoint. 
Many functors throughout algebra turn out to be examples of such left 
adjoints. Other highlights of functorial semantics include the reconstruction 
of the theory $T$ from its category of algebras in sets, $Alg(T, Sets)$. We 
also learn how to combine two algebraic structures using the Kronecker product 
construction \cite{Freyd}.

This paper deals with the two-dimensional analog of functorial 
semantics. We start with the definition of an algebraic 2-theory, $\TT$, 
\cite{Grayb} which is a 2-category whose objects are the natural numbers, 
whose morphisms correspond to operations (functors) and whose 
2-cells correspond to 
natural transformations between functors. We then go on to define a 
2-theory-morphisms and other morphisms in $\TTh$ (following \cite{Gray},
we place a tilde over all 3-categories.)  Connections between $\Th$ and $\TTh$ 
are enumerated.

Algebras for $\TT$ are product preserving functors from $\TT^{op}$ to a 
2-category, $\C$, with a product structure. If $\C$ is $\Cat$ then algebras are 
categories with extra structure. Not all structures that are put on a category 
can be represented by a 2-theory. We are restricted to structures with only 
covariant functors and hence can not deal with a closure or a duality 
structure. Methods of generalizing this work in order to handle such 
structures are discussed in section 5.

A search through the literature reveals that morphisms between algebras generally 
do not preserve the operations ``on the nose.'' Rather, They are preserved up to a natural 
(iso)morphism. This translates to the notion of a quasi-natural 
transformation \cite{Bunge1, Bunge2, Gray}
between the product preserving functors from $\TT^{op}$ to $\C$
i.e. a naturality square that commutes up to a 2-cell. It is important to realize that 
the quasi-natural transformations places our subject outside of
enriched functorial semantics \cite{Bor&Day}.  
Between quasi-natural transformations there are modifications/2-cells. And
so we have the 2-category ${\bf 2Alg(T,C)}$ of algebras, quasi-natural transformations
and modifications.

Section 3 discusses the left quasi-adjoint $Lan_G$ of $G^*: {\bf 
2Alg(T',C)} \longrightarrow {\bf 2Alg(T,C)}$ where $G:\TT \longrightarrow \TT'$ is a 
2-theory-morphism.
In order to construct this quasi-left-Kan extension we must first talk of 
the quasi-Yoneda lemma, quasi-comma categories, quasi-cocones, quasi-colimits etc.
Our aim is not 
to repeat all the superb work of \cite{Bunge1, Bunge2, Gray, Street, Bourn, Mac&Sto} on quasi-(co)limits
and quasi-adjoints, rather it is to state only what we need for 
functorial semantics. We have aimed at making this as readable as possible and 
we do not assume knowledge of any of the above papers. This paper is self-contained.
The main idea behind section 3 is given 
two 2-theories and a 2-theory-morphism between 
them, one should explore the induced relationship 
between the corresponding 2-categories of 
algebras. The strength of the induced 
quasi-adjoints are classified by the strength of 
the 2-theory-morphism. These quasi-adjoints 
reflect the extent to which one structure can be 
replaced by another. Different types of 2-theory-morphisms
induce  quasi-adjoints of varying strength and these different adjoints 
express the coherence results. Whereas in the 1-dimensional case, if the left adjoint is
an equivalence of categories, the theories are isomorphic, in the 2-dimensional case, there are 
many intermediate possibilities. The aim is to simply look at the combinatorics 
of the 2-theory-morphism in order to understand the coherence result that is implied.
Many examples are given. We also show how to reconstruct the 2-theory from the 
2-category of algebras.

Section 4 is a discussion of a two-dimensional generalization of the Kronecker 
product. We show how one can combine one structure  with another. This leads 
to many examples that are the standard fare of coherence theory. We go on to 
see how this Kronecker product respects the left quasi-adjoints of section 3. 
This helps us combine coherence results.

We end the paper with a look at the different directions that this project 
can proceed. Several conjectures are made. Some questions that seem interesting 
and important for future work are asked.  
Applications to representation theory and
physics are discussed.

This paper was written to be self-contained. We assume only the basic 
definitions of 2-category theory. However, this work does not stand alone. 
This paper --- as all papers in higher category theory --- owes much to John 
Gray's important ground-breaking book \cite{Gray}. We try to follow his names 
and notation when possible. Many of our examples come from Joyal and Street's 
wonderful paper \cite{Joy&Str} on the many structures and coherence theorems 
that are important for modern mathematics. This project would not exist 
without either of these important works. 

A note on notation. In order to alleviate the pain of all the different 
types of morphisms, we will call objects and morphisms ``0-cells'' and ``1-cells'' 
respectively. However, an $i$-cell in one category can be an $i'$-cell in another 
category. All 2-categories will be in bold typeface. In contrast, (1-)categories 
will not.

This work would not have been possible without many enlightening and 
encouraging conversations with Alex Heller. Some of the ideas for this paper 
were formed while enjoying a postdoctoral position at McGill University. I 
would like to thank Mike Barr, Marta Bunge, Michael Makkai, Robert Seely
and the entire Montreal categories group for that wonderful experience. 
Jim Stasheff and Mirco Mannucci have looked at earlier drafts of this 
paper and made many helpful comments. I am indebted to them.

\section{2-theories and 2-algebras}
Consider the skeletal category of finite sets $\fin$. Place a
coproduct structure on this category. The coproduct structure
allows one to write $n \cong \coprod_n 1$. For all positive
integers $m, n$ and $p$, one has the following induced
isomorphisms
$$\sigma_m^n: \coprod_m n \cong \coprod_m \coprod_n 1
\longrightarrow \coprod_n \coprod_m 1 \cong \coprod_n m$$
$$ \mu_{m,n}^p: \coprod_m p + \coprod_n p \longrightarrow
\coprod_{m+n} p$$
$$\nu_p^{m,n}: \coprod_p m + \coprod_p n \longrightarrow
\coprod_p (m+n).$$
These isomorphisms satisfy the following coherence condition:
\be 
\label{XY:coprod}
\xymatrix{
\coprod_m p + \coprod_n p \ar[rr]^{\mu^p_{m,n}}
\ar[dd]_{\sigma + \sigma}
 &&
\coprod_{m+n} p
\ar[dd]^{\sigma}
\\
\\
\coprod_p m + \coprod_p n
\ar[rr]^{\nu_p^{m,n}}
&&
\coprod_p {m+n}.}
\ee

Let $\finbar$ denote the 2-category with the same 0-cells and 
1-cells as $\fin$ but with only identity 2-cells. $\finbar$ 
also has a coproduct structure. A coproduct structure for a 
2-category is similar to a coproduct structure for a 1-category. However, 
there is an added requirement that for every finite family of
1-cells with common source and target, there is a 1-cell with 
injection 2-cells that satisfy the obvious universal property.
When we talk of preserving coproduct structures, we mean preserving
the coproduct strictly (equality). 
 
\begin{defi}
A (single sorted algebraic) 2-theory is a 2-category $\TT$  
 with a given coproduct structure and a 2-functor $G_{\TT}:
\finbar \longrightarrow \TT$ such that $G_{\TT}$ is bijective on 
0-cells and preserves the coproduct structure.
\end{defi}

The following examples are well known.

\begin{examp}
$\finbar$ is the initial 2-theory. Just as $\fin$ is
the theory of sets, so too,  $\finbar$ is the theory of
categories.
\end{examp}

\begin{examp}
Let $\binbar$ be $\finbar$ with a nontrivial generating 
1-cell $\otimes:1 \longrightarrow 2$ thought of as a binary
operation (bifunctor).
\end{examp}

\begin{examp}
$\Tmon$ is the 2-theory of monoidal (tensor) categories. It is a 2-theory 
``over'' $\binbar$ with a 1-cell $e:1 \longrightarrow 0$.
The isomorphic 2-cells are generated by 
$$
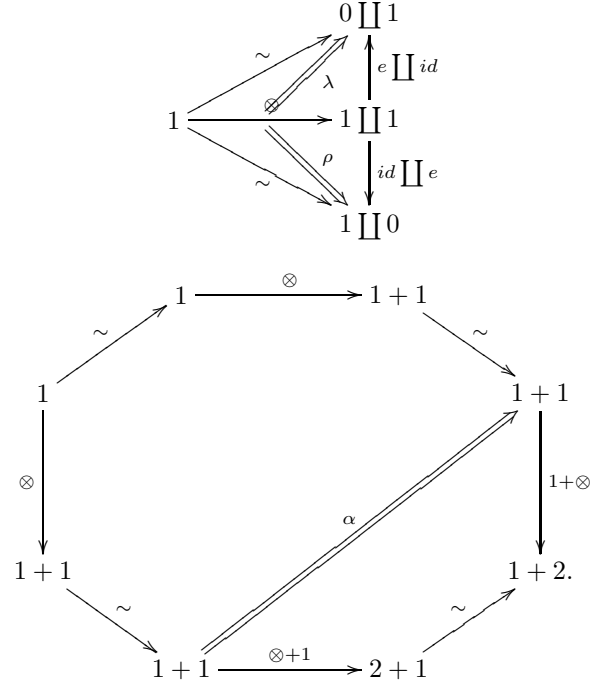
\xymatrix{
&& 0 \coprod 1
\\
1 \ar[urr]^{\sim} \ar[rr]^{\otimes} \ar[drr]_{\sim}
&{ }\ar@{=>}[ur]_{\lambda}\ar@{=>}[dr]^{\rho}&
 1 \coprod 1 \ar[u]_{e \coprod id}  \ar[d]^{id \coprod e}
\\
&& 1 \coprod 0
}$$
$$\xymatrix{
     &1\ar[rr]^{\otimes} & & 1+1 \ar[dr]^{\sim}
\\
1\ar[ur]^{\sim}\ar[dd]_{\otimes}&  & & &1+1 \ar[dd]^{1+\otimes}
\\
\\
1+1\ar[dr]^{\sim} &&&& 1+2.
\\
& 1+1 \ar[rr]^{\otimes + 1}\ar@{=>}[uuurrr]^{\alpha} && 2+1 \ar[ur]^{\sim}
}$$
where the corner isomorphisms $n+m \lra m+n$ is an instance 
of $\sigma^n_m$ in $\finbar$. 
These 2-cells are subject to a unital equation (left for the
reader) and the now-famous pentagon condition:
$$
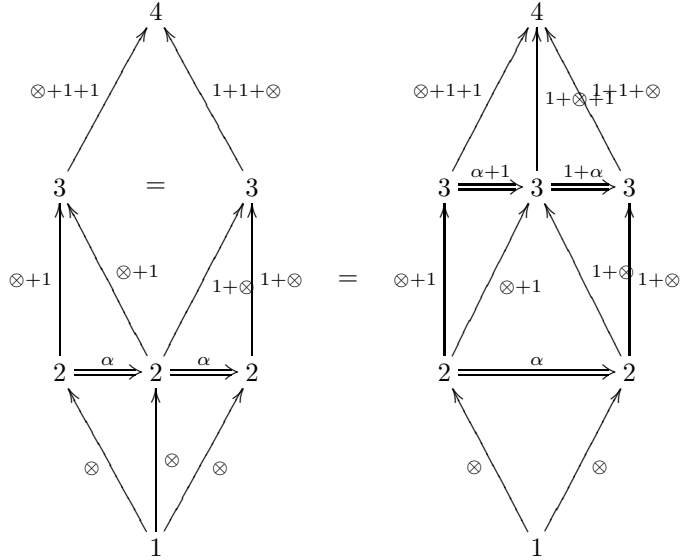
\xymatrix{
      &   4   &      &        &      &   4
\\
\\
3\ar[uur]^{\otimes +1+1}&=       &3\ar[uul]_{1+1+\otimes}  &        &3\ar@{=>}[r]^{\alpha+1}\ar[ruu]^{\otimes+1+1}  &
3\ar@{=>}[r]^{1+\alpha}\ar[uu]_{1+\otimes+1} &3\ar[uul]_{1+1+\otimes}       
\\    &       &      &   =    &
\\
2\ar[uu]^{\otimes+1}\ar@{=>}[r]^{\alpha}&2\ar@{=>}[r]^{\alpha}\ar[uur]_{1+\otimes}\ar[uul]_{\otimes +1}  &2\ar[uu]_{1+\otimes}
  &&2\ar@{=>}[rr]^{\alpha}\ar[uu]^{\otimes+1}\ar[uur]_{\otimes+1}  &       & 2\ar[uul]_{1+\otimes}
\ar[uu]_{1+\otimes}
\\
\\
      &   1\ar[luu]^{\otimes} \ar[uu]_{\otimes} \ar[ruu]_{\otimes}   &      &     
   &      &   1 \ar[luu]^{\otimes}
 \ar[ruu]_{\otimes}  
}$$
(We leave out the corner isomorphisms in order to make the
diagram easier to read. However they are important and must be
placed in the definition).
\end{examp}

\begin{examp}
The theory of braided tensor categories $\Tbraid$ and balanced
tensor categories \cite{Joy&Str},$\Tbal$,  are easily 
described in a similar manner.
\end{examp}

\begin{examp}
Associative categories \cite{Paper1} which are monoidal categories 
in which the pentagon coherence does not necessarily hold are 
described by $\Tassoc$. Similarly, commutative categories \cite{Paper2}
which are braided tensor categories that do not necessarily satisfy the
hexagon coherence condition are described by $\Tcomm$.
\end{examp}

\begin{examp}
Whenever we have a theory with strict associativity, we denote it with
a small ``s'' followed by the usual name e.g. $\Tsmon$, $\Tsbraid$,
$\Tsbal$ etc.
\end{examp}

\begin{defi}
A 2-theory-morphism from $\T 1 $ to $\T 2 $ is a 2-functor $G:\T 1
\longrightarrow \T2 $ such that
$$
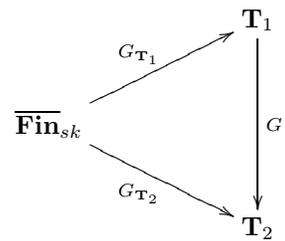
\xymatrix{
&& \T1 \ar[dd]^G
\\
\finbar \ar[urr]^{G_{\T 1 }}
\ar[drr]_{G_{\T 2 }}
\\
&& \T 2
}$$
commutes.
\end{defi}

\begin{defi} A 2-theory-natural transformation $\gamma :G_1
\Longrightarrow G_2$ between two 2-theory-morphisms is a natural
transformation such that 
$$
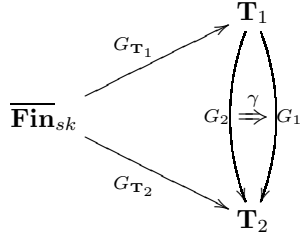
\xymatrix{
&& \T 1 \ddtwocell^{G_1}_{G_2}{^\gamma}
\\
\finbar \ar[urr]^{G_{\T 1 }}
\ar[drr]_{G_{\T 2 }}&
\\
&& \T 2
}$$
commutes.
\end{defi}

One goes on to define a 2-theory-modification $\Gamma :\gamma_1
\leadsto \gamma_2$ in the obvious way.

We shall denote the 3-category of 2-theories, 2-theory-morphisms,
natural transformations and modifications as $\TTh$.

Here is a diagram of some of the 2-theories and 2-theory 
morphisms that we will work with.
$$ 
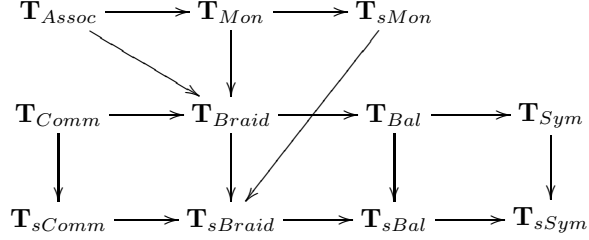
\xymatrix{
\Tassoc \ar[r] \ar[dr] & \Tmon \ar[r] \ar[d] & \Tsmon \ar[ddl]
\\
\Tcomm \ar[r] \ar[d]& \Tbraid \ar[r]\ar[d] & \Tbal \ar[r]\ar[d] & \Tsym \ar[d]
\\
\Tscomm \ar[r] & \Tsbraid \ar[r] & \Tsbal \ar[r] & \Tssym
}$$

Many examples of 2-theories  and their morphisms come from one 
dimensional theories in the following manner. Let $\th $ denote
the usual \cite{Lawvere} 2-category of theories,
theory-morphisms and theory-natural transformations.
One can think of $\th $ as a 3-category 
$\Th $ with only trivial 3-cells. Analogous to the
relationship between sets and topological spaces, we have the
following adjunctions:
$$
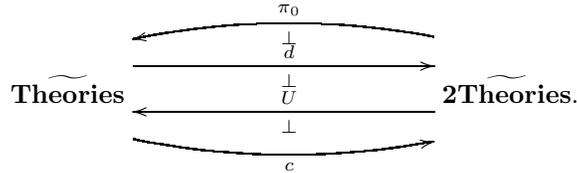
\xymatrix{
\Th \ar@/_1pc/@<-3ex>[rrrr]_c 
\ar@<2ex>[rrrr]^d_{\bot}
&&&&\TTh .
\ar@/_1pc/@<-3ex>[llll]_{\pi_0}^{\bot}
\ar@<2ex>[llll]_U^{\bot}
}$$

$c(T)$ is the 2-theory with the same 1-cells as $T$ and a unique
2-cell between nontrivial 1-cells. $d(T)$ has the same 1-cells as
$T$ and only trivial 2-cells. $U( \TT)$ forgets the 2-cells of
$\TT$. $\pi_0( \TT)$ is a quotient theory of $\TT$ where two 1-cells
are set equal if there is a 2-cell between them. These functors
extend in an obvious way to 3-functors. By adjunction we mean a
strict 3-adjunction; that is the universal property is satisfied
by a strict 2-category isomorphism. For example the following 2-categories
 are isomorphic
$$Hom_{\Th}(T, U(\TT)) \cong Hom_{\TTh}(d(T), \TT)$$

\begin{examp}
$\finbar = d(\fin )$, that is, the theory of categories is
the discrete theory of sets.
\end{examp}

\begin{examp}
$\binbar = d(T_{Magmas})$.
\end{examp}

\begin{examp}
$d(T_{Monoids})$ is the theory of strict monoidal
categories, $\Tsmon$. 
\end{examp}

\begin{examp}
Let $T_{Magmas \bullet}$ be the theory of pointed magmas
i.e. the theory of magmas with a distinguished element.
$c(T_{Magmas \bullet})$ is the 2-theory of symmetric
(monoidal ) tensor categories. Warning: not all operations
are made to be isomorphic. In particular, the projections (inclusions)
live in $\finbar$ and are not isomorphic.
\end{examp}

\begin{examp}Let $\Tbraid $ denote the 2-theory of braided tensor categories.
$\pi_0(\Tbraid )$ is
the theory of commutative monoids.
\end{examp}

The units and counits of these adjunctions are of interest. 
$\varepsilon : \pi_0 d T \longrightarrow T$,
$\mu : T \longrightarrow U d T$ and 
 $ \varepsilon : U c T \longrightarrow T$ are all identity theory-morphisms.
More importantly, $\mu : \TT \longrightarrow d \pi_0 \TT$ is the 
2-theory-morphism corresponding to ``strictification''. Every 2-cell becomes the
identity. ``Strictification'' is often used in coherence theory and in section 3 we shall
take (quasi-) Kan extensions along such 2-theory-morphisms. Similarly, 
 $\mu : \TT \longrightarrow c U \TT$
might be called ``coherification'': a 2-theory is forced to be coherent. 
$\varepsilon: d U \TT \longrightarrow \TT$ is the injection of the 1-theory into
the 2-theory.

Given a 2-theory $\TT$ and a 2-category $\C $ with a product
structure, an {\bf algebra} of $\TT$ in $\C $ is a product preserving 2-functor
 $F:\TT^{op} \longrightarrow \C. $

A {\bf quasi-natural transformation} (cf. pg. 26 of \cite{Gray}, \cite{Bunge1, Bunge2})
 $\sigma$ from an algebra $F$ to an to an algebra $F'$ is 
\begin{itemize}
\item A family of 1-cells in $\C $, $\sigma_n:F(n) \longrightarrow
F'(n)$ indexed by 0-cells of $\TT $. This family must preserve
products i.e. $\sigma_n = (\sigma_1)^n:F(1)^n \longrightarrow
F'(1)^n$.
\item A family of 2-cells in $\C $, $\sigma_f$, indexed by 1-cells
$f:m\longrightarrow n$ of $\TT$. $\sigma_f$ makes the following diagram commute.
\end{itemize}
\be 
\label{XY:quasi}
\xymatrix{
& F(1)^n \ar[r]^{\sigma^n}& F'(1)^n\ar[dr]^{\sim} 
\\
F(n)\ar[ur]^{\sim} \ar[dd]_{Ff}& & & F'(n)\ar[dd]_{F'f}&
\\
&&& { }
\\
F(m) \ar[dr]_{\sim} &&& F'(m) &
\\
& F(1)^m \ar[r]^{\sigma^m}\ar@{=>}[rruu]^{\sigma_f}  & F'(1)^m \ar[ur]_{\sim} &&
}\ee

These morphisms must satisfy the following conditions:
\begin{enumerate}
\item If $f$ is in the image of $G_\TT:\finbar \longrightarrow \TT$,
then $\sigma_f = id$. That is, diagram (\ref{XY:quasi}) commutes strictly. This
condition includes $\sigma_{id_n} = id_{\sigma_n}$.
\item $\sigma$ preserves the coproduct structure: $\sigma_{f + f'} =
 \sigma_f \times \sigma_{f'}$.  To be more exact, $\sigma_{f + f'}$ is the 
entire diagram in Figure I.
The quadrilaterals in Figure I commute from the coproduct structure of $\TT$ 
and the product structure of $\C$; see diagram  (\ref{XY:coprod}).
\item $\sigma_{g\circ f} = \sigma_f \circ_v \sigma_g $where $\circ_v$ is the 
vertical composition of 2-cells.
\item $\sigma$ behaves well with respect to 2-cells of $\TT$. That is, if we 
have 
$$\xymatrix{
m \rrtwocell^f_{f'}{\alpha}&&n 
} $$
in $\TT$, then the two diagrams of Figure II must be equal.
\end{enumerate}
\vspace{1cm}

\tiny
$$
\xymatrix{
&F(1)^{n+n'} \ar[rrrr] \ar[dr] &&&& F'(1)^{n+n'} \ar[dr]
\\
F(n + n') \ar[dr] \ar[ru] \ar[dddd]_{F(f+f')}&&F(1)^n\times F(1)^{n'} 
 \ar[rr]^{(\sigma_1)^n\times(\sigma_1)^{n'}  } & &F'(1)^n \times  F'(1)^{n'}
\ar[dr] \ar[ur] & & F'(n+n') \ar[dddd]_{F'(f+f')}
\\
&F(n)\times F(n') \ar[ur]^{ }\ar[dd]_{F(f)\times F(f')}&  & & &F'(n)
\times F'(n') \ar[dd]_{F'(f)\times F'(f')} \ar[ur]
\\
\\
&F(m)\times F(m') \ar[dr]^{ }&  & & &F'(m)\times F'(m') \ar[dr]           
\\
F(m + m')\ar[ur] \ar[dr]&& F(1)^m\times F(1)^{m'}\ar@{=>}[rrruuu]^
{\sigma_f\times \sigma_{f'}}
 \ar[rr]^{(\sigma_1)^m\times(\sigma_1)^{m'}} 
&& F'(1)^m\times F'(1)^{m'} \ar[dr]
\ar[ur]^{ }&&F'(m+m')
\\
&F(1)^{m+m'} \ar[rrrr] \ar[ru] &&&& F'(1)^{m+m'} \ar[ur]
}$$
\normalsize
\centerline{\bf Figure I.}
\vspace{2cm}
\tiny
$$
\xymatrix{
& F(1)^n \ar[r]^{\sigma^n}& F'(1)^n\ar[dr] && & F(1)^n \ar[r]^{\sigma^n} & F'(1)^n\ar[dr]
\\
F(n)\ar[ur] \ddtwocell_{Ff}^{Ff'}{^F\alpha}& & & F'(n)\ar[dd]_{F'f'}&  F(n)\ar[dd]_{Ff}\ar[ur] & &
 & F'(n) \ddtwocell_{F'f}^{F'f'}{^F'\alpha}
\\
&&&\ar@{}[r]^=& &  & &
\\
F(m) \ar[dr] &&& F'(m) & F(m)\ar[dr] & & & F'(m)
\\
& F(1)^m \ar[r]^{\sigma^m}\ar@{=>}[rruu]^{\sigma_f}  & F'(1)^m \ar[ur] &&
  & F(1)^m\ar@{=>}[rruu]^{\sigma_{f'}}  \ar[r]^{\sigma^m} & F'(1)^m. \ar[ur]}$$
\normalsize 
\centerline{\bf Figure II.}

\begin{remar}We not only require $\sigma$ to preserve the coproduct
in $\TT$ but also to preserve {\em all} the coherence properties
of the coproduct.
\end{remar}

Composition of quasi-natural transformations are  given as 
$$ (\sigma' \sigma)_n = \sigma'_n \sigma_n \qquad (\sigma'\sigma)_f = \sigma'_f 
\circ_h \sigma_f$$

Given two quasi-natural transformations $\sigma , \sigma' : F \longrightarrow 
F'$, a
{\bf modification} $\Sigma:\sigma \leadsto \sigma'$ from $\sigma$ to $\sigma'$
is a family of 2-cells $\Sigma_n:\sigma_n \Longrightarrow \sigma'_n$ indexed by
the 0-cells of $\TT$. These 2-cells must satisfy the following conditions:
\begin{enumerate}
\item $\Sigma$ preserves products i.e. $\Sigma_n = (\Sigma_1)^n:(\sigma_1)^n 
\Longrightarrow (\sigma_1')^n$.
\item $\Sigma$ behaves well with respect to the 2-cells of $\TT $. That is,
if we have 
$$\xymatrix{
m\rrtwocell^f_{f'}{\alpha}&&n
} $$
then we have the following ``cube relation'':
\end{enumerate}
\be
\label{XY:qcube}
\xymatrix{
F(n) \ar[rr]^{\sigma'_n} \ar[dd]_{id} && F'(n) \ar[dr]^{F'(f')} \ar[dd]^{id}
 &&  
F(n) \ar[rr]^{\sigma'_n} \ar[dd]_{id} \ar[dr]^{F(f')} && F'(n) \ar[dr]^{F'(f')} 
\\
 & &{ } & F'(m) \ar[dd]^{id} 
& &        F(m) \ar[rr]_{\sigma'_m} \ar[dd]^{id} &{ }\ar@{=>}[u]_{\sigma'_{f'}} & F'(m) \ar[dd]^{id}
\\
F(n) \ar[rr]_{\sigma_n} \ar[rd]_{F(f)} &{ }\ar@{=>}[ur]^{\Sigma_n}& F'(n) \ar[rd]^{F'(f)}\ar@{=>}[r]^{F'(\alpha)}
&{ }\ar@{}[r]^= & F(n) \ar[dr]_{F(f)}\ar@{=>}[r]^{F(\alpha)}&&&{ }
\\
& F(m) \ar[rr]_{\sigma_m} &{ }\ar@{=>}[u]^{\sigma_f} &F'(m)
&& F(m) \ar[rr]_{\sigma_m} &{ }\ar@{=>}[ur]^{\Sigma_m}&F'(m)
}\ee

Compositions of modifications are given as 
$$ (\Sigma' \circ_h \Sigma)_n = \Sigma'_n \circ_h \Sigma_n \qquad
 (\Sigma' \circ_v \Sigma)_n = \Sigma'_n \circ_v \Sigma_n $$

There is a need to generalize this definition. Let  $G:\T 1 \longrightarrow \T 2$
be a 2-theory-morphism. Then ${\bf 2Alg}_G(\T 2,\C )$ will have the same
 0-cells as ${\bf 2Alg(\T 2,\C )}$, however,   ${\bf 2Alg_G(\T 2,\C )}(F, F')$ will be the 
full subcategory of  ${\bf 2Alg(\T 2,\C)}(F,F') $ consisting of those 
quasi-natural transformations that are actual natural transformations when 
precomposed with $G$ i.e. those $\sigma$ such that $\sigma_{G(f)} = id$ or
in other words those $\sigma$ such that 
$$ \xymatrix{
\T 1^{op}\ar[rr]^{G^{op}} && \T 2^{op}\rrtwocell^F_{F'}{\sigma} && \C }$$
 $\sigma \circ G^{op}$ is a natural transformation (not quasi)  from $F \circ G^{op}$  
to  $F' \circ G^{op}$.

$\TA^i_G(\T 1,\C)$ is defined to be the locally full sub-2-category of  
$\TA_G(\T 1,\C)$ consisting of those quasi-natural transformations
where the $\sigma_f$'s  are  isomorphisms.

It is obvious that ${\bf 2Alg_{G_\TT}(\TT,\C)=2Alg(\TT,\C)}$ and that    
$\TA_{id_\TT}(\TT,\C)$ has only strict natural transformations.
For every ${\bf 2Alg_G(\TT,\C)}$, there is a forgetful 2-functor 
 $U:{\bf 2Alg_G(\TT,\C)} \longrightarrow \C$ defined as $U(F) = F(1)$;
 $U(\sigma)=\sigma_1$ and $U(\Sigma)=\Sigma_1$.

Consider the 3-category $(\TTh)^{\rightarrow}$ which has as 0-cells 
2-theory-morphisms. The $i$-cells for $i=1,2,3$ are pairs of $i$-cells in 
$\TTh$ making the usual square commute. Thus we have the following 3-functor
$$\TA_{(?)}(t(?), \C):((\TTh)^{\rightarrow})^{op} \longrightarrow {\bf 
2Cat/\C}$$ where $t(?)$ is the target (codomain) of $(?)$.

\section{Universal properties of coherence}

Many coherence theorems are a result of the quasi-Yoneda lemma.
\begin{lem}[Quasi-Yoneda] Let $\D$ be a 2-category. Let $K:\D \longrightarrow \Cat$
be a 2-functor. ${\bf qNat} ( \D (r,-), K(-))$ shall denote the category of 
quasi-natural transformations and modifications (not necessarily product preserving) 
between $\D (r,-)$ and $K(-)$. Then there are (quasi-)adjoint functors:
$$  
\xymatrix{ 
{\bf qNat}(\D (r,-), K(-)) \ar@<1ex>[rr]^\Psir_\bot && \qquad K(r). \ar@<1ex>[ll]^\Psirh
}$$
The unit of this adjunction, $id \longrightarrow \Psir \circ \Psirh$, is quasi-natural.
The counit of this adjunction, $\Psirh \circ \Psir \longrightarrow id$, is the identity.
\end{lem}
{\bf Proof.} \underline{ Definition of $\Psir$.} Let $\sigma: \D (r,-) \longrightarrow K(-)$
be a quasi-natural transformation then $\Psir(\sigma)=\sigma_{r,id_r} \in K(r)$. For
a modification $\Sigma: \sigma  \leadsto \sigma'$, we set $\Psir(\Sigma)=\Sigma_{r,id_r}$
where $\Sigma_r$ is a 2-cell in $\Cat$ (a natural transformation):
$$
\xymatrix{
\D (r,r) \rrtwocell^{\sigma_r}_{\sigma'_r}{\Sigma_r} && K(r). 
}
$$

\underline{ Definition of $\Psirh$.} Let $U \in K(r)$. $\Psirh(U)= \sigma_U$ where 
$\sigma_{U,d}:\D (r,d) \longrightarrow K(d)$ is defined as follows. For $f \in
\D(r,d)$, $\sigma_{U,d}(f)=K(f)U \in K(d)$. For $\alpha: f \Longrightarrow f'$,
 $\sigma_{U,d}(\alpha)=K(\alpha)U$. One should have the following picture in mind:
$$\xymatrix{
\underline{\D} &&& \underline{\Cat}
\\
r \ddtwocell_f^{f'}{^\alpha} &&& {U \in K(r)} \ddtwocell_{Kf}^{Kf'}{^K\alpha}
\\
&\ar@{|~>}[r]& &{ }
\\
d &&& K(f)U \in K(d) \ni K(f')U.}$$
Let $t:U\longrightarrow U'$ be a 1-cell in $K(r)$. Then $\Psirh(t)=\Sigma_t$
where $\Sigma_{t,d}$ fits in 
  $$
\xymatrix{
\D (r,d) \rrtwocell^{\sigma_U}_{\sigma_{U'}}{\Sigma_{td}} && K(d). 
}
$$ 
and is defined as follows:
$$[\Sigma_{t,d}(f: r \rightarrow d) = K(f)(t)]\quad :\quad [\sigma_{U,d}(f)=K(f)U]
 \quad \longrightarrow 
\quad [\sigma_{U',d}(f)=K(f)U'].$$
And finally
$$\Sigma_{t,d}(\alpha:f \Rightarrow f') \quad = \quad K(f')t \circ K(\alpha)_U \quad =
 \quad K(\alpha)_{U'} \circ K(f)t$$
i.e. the morphism described by the natural transformation
$$\xymatrix{
K(f)U \ar[rr]^{K(\alpha)_U}  \ar[dd]_{K(f)t} && K(f')U \ar[dd]^{K(f')t}
\\
\\
K(f)U' \ar[rr]^{K(\alpha)_{U'}} \ar[rr] && K(f')U'.
}$$

\underline{The unit of the adjunction.} Let $\sigma$ be a quasi-natural transformation.
$\Psir(\sigma)=\sigma_{r,id_r}.$ Then $\Psirh \Psir(\sigma) = \Psirh(\sigma_{r,id_r})$.
$$ \Psirh(\sigma_{r,id_r})_d : \D(r,d) \longrightarrow K(d)$$
 is defined as 
$$ \Psirh(\sigma_{r,id_r})_d (f) = K(f)(\sigma_{r,id_r})$$ The unit of the adjunction
at $d\in \D$, $\sigma_d \longrightarrow \Psirh \Psir(\sigma)_d$ , is defined at 
$f\in \D (r,d)$ as $\sigma_{f,id_r}$. The following picture is helpful:
$$
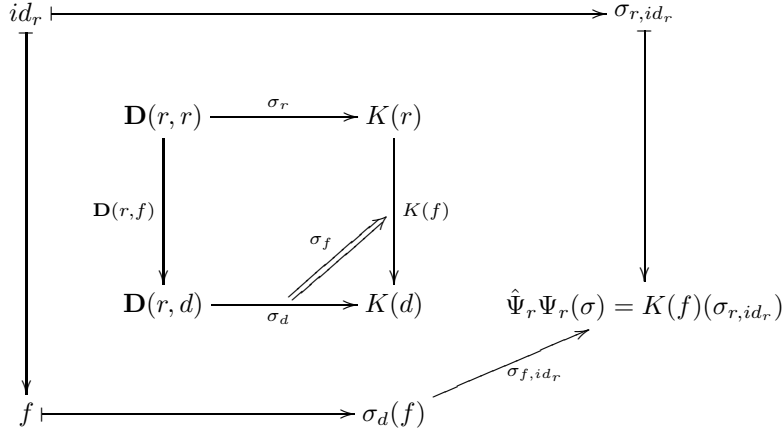
\xymatrix{
id_r \ar@{|->}[rrrr] \ar@{|->}[dddd] &&&& \sigma_{r,id_r} \ar@{|->}[ddd]
\\
& \D (r,r) \ar[rr]^{\sigma_r} \ar[dd]_{\D (r,f)} && K(r) \ar[dd]^{K(f)}
\\
 & & & { }
\\
&\D (r,d) \ar[rr]_{\sigma_d} & { }\ar@{=>}[ru]^{\sigma_f} & K(d) &  \Psirh \Psir(\sigma) = K(f)(\sigma_{r,id_r})
\\
f \ar@{|->}[rrr] &&& \sigma_d(f) \ar[ru]_{\sigma_{f,id_r}}}
$$

Note that if we insist that $\sigma_f$ is an isomorphism, then $\sigma_{f,id_r}$
is also an isomorphism and hence the unit would be an isomorphism. The unit 
is a quasi-natural transformation.

\underline{The counit of the adjunction.} 
$$ (\Psir \circ \Psirh)(U) \quad = \quad \Psir(\sigma_U) \quad = \quad \sigma_{U,r,id_r} 
\quad = \quad K(id_r)U \quad$$
$$ = \quad id_{K(r)}(U) \quad = \quad U \qquad \Box $$

This theorem says that every 0-cell in $K(r)$ corresponds
to a natural transformation. The unit of the adjunction is,
in a sense, a reflection of the category of natural transformations inside 
the category of quasi-natural transformations.

The following facts about the quasi-Yoneda lemma are important. The 
proofs are trivial or tedious and we leave them for the 
readers leisure time.
  
\begin{prop}[On the quasi-Yoneda lemma] Let $l:r \longrightarrow r'$
be a 1-cell in $\D$ and let $\kappa:K \longrightarrow K'$ be a quasi-natural 
transformation.
\renewcommand{\labelenumi}{(\alph{enumi})}

\begin{enumerate}
\item $\Psir$ is quasi-natural with respect to $ r $
i.e. given  $l : r \longrightarrow r'$, the obvious square commutes up
to a natural transformation.
\item $\Psir$ is natural with respect to a $K$ i.e. given a quasi-natural transformation
$\kappa:K \longrightarrow K'$, the obvious square commutes strictly.
\item $\Psirh$ is natural with respect to $r$.
\item $\Psirh$ is quasi-natural with respect to $K$.
If however, $\kappa$ is natural (not quasi), then $\Psirh$ is also natural
(not quasi).
\item If we insist that the quasi-natural transformations $\sigma$ have the 
usual square commuting up to a natural {\bf iso}-2-cell, then the $\Psir, \Psirh$ 
adjunction becomes an equivalence of categories:
$${\bf qNat^i} (  \D (r,-) , K(-)) \cong K(r). $$
Warning: This is not natural in $r$.
\item If we insist that the quasi-natural transformations $\sigma$ be   
$\Cat$-natural transformations, then the $\Psir, \Psirh$ 
adjunction becomes an isomorphism of categories:
$${\bf CatNat} (  \D (r,-) , K(-)) \cong K(r). $$
\item If we insist that all the 2-cells are identities, then the $\Psir, \Psirh$ 
adjunction become the usual Yoneda lemma:
$${Nat} (  D (r,-) , K(-)) \cong K(r). $$
\item If $\D$ has a product structure and $K, \sigma, \Sigma$ are assumed to preserve 
the product structure, then we still have the adjunction. Furthermore,
$ \Psi_{r\times r'} \cong \Psir \times \Psi_{r'} $  (similarly for $\Psirh$).
\end{enumerate}
\end{prop}

Let
$[n]$ denote the discrete category whose objects are $\{0, \ldots, n-1\}$.
\begin{prop} Let $\TT$ be a 2-theory. $\TTo(n,-):\TTo 
\longrightarrow \Cat$ is a product preserving 2-functor
and is the free $\TT$-algebra on $n$ generators in the 
sense that 
$$\TAi(\TT,\Cat)(\TTo(n,-),F(-)) = {\bf qNat^i}((\TTo(n,-),F(-)) \cong
F(n)$$
$$  \cong F(1)^n \cong \Cat([n], F(1)) \Box . $$
\end{prop}
Notice the importance of insisting on {\it iso}-quasi-natural transformations
since by (g) above, we have an equivalence of categories.
From the universality and (quasi-)naturality of the quasi-Yoneda lemma, 
any other $\TT$-algebra that satisfies this universal property is
equivalent to $\TTo(n,-)$ in $\Cat$  and is equivalent to $\TTo(n,-)$ in
$\TAi(\TT,\Cat)$.

\begin{examp} Let $B$ (see page 10 of \cite{Joy&Str0}) 
be the category whose objects are the natural numbers
and whose only morphisms are $Hom_B(n,n)=B_n$, 
the Artin braid group on $n$ strings. $B$ has a
strict braided structure and is the free strict braided
tensor category on one generator. Let $\Tsbraid$ be the 2-theory of strict
braided tensor categories.
$\Tsbraid^{op}(n,-)$ is the free braided tensor category generated by $n$
objects Hence $B \cong \Tsbraid^{op}(1,-)$
as categories and as braided tensor categories.
\end{examp}

\begin{examp}
Let $S$ be the category similar to $B$ but whose morphisms are the 
symmetric groups. $S$ has a strict symmetric structure and is the 
free strict symmetric tensor category on one generator.
Let $\Tssym$ be the 2-theory of strict symmetric tensor categories. Hence $S \cong 
\Tssym^{op}(1,-)$ as categories and as symmetric tensor categories.
\end{examp}

\begin{examp}
Let $\tilde{B}$ be the 
free strict balanced  tensor category on one generator (see pages 11, 41
of \cite{Joy&Str0}). Let $\Tbal$ be the 2-theory of strict balanced tensor 
categories. Hence $\tilde{B} \cong \Tbal^{op}(1,-)$. 
\end{examp}

On to the notion of quasi-cocones. Let $\I$ be a small, locally small 2-category. 
$\CI$ shall denote the 2-category 
of $\I$-diagrams in $\C$. The 1-cells in $\CI$ are quasi-natural transformations.
In order to keep track of the morass of different
types of morphisms in this discussion, we shall attempt to abide by the following table.
$$
\begin{tabular}{l||c c c} 
			& $\I$		& $\C$& $\CI$ \\
\hline \hline
\mbox{0-cells}		& $i$		& $c$	& $d$ \\ \hline
\mbox{1-cells}	 	&  $I  :  i  \longrightarrow i'$	& $\sigma : c\longrightarrow c'$	&$\xi : d \longrightarrow
d'$\\ \hline
\mbox{2-cells}		&$\iota : I \Longrightarrow I'$	& $\Sigma : \sigma \Longrightarrow \sigma'$ & 
 $\Xi : \xi \Longrightarrow \xi' $\\ \hline
\end{tabular}
$$
\\

To every $\I$, there is a constant-diagram 2-functor
$$\Delta : \C \longrightarrow \CI$$
which is defined on 0-cells as follows
$$\Delta(c)(i) = c \qquad \Delta(c)(I) =id_c \qquad \Delta(c)(\iota) = id_{id_c}.$$
$\Delta(\sigma)$ and $\Delta(\Sigma)$ are defined to be the usual morphisms 
between constant 2-diagrams.

The category $\CI(d, \Delta(c))$ is the category of cocones over $d$ with vertex $c$ and
 morphisms between such cocones. In detail, a cocone $\xi$ over $d$ with vertex $c$ is a 
quasi-natural transformation in the 2-category  $\C$. For every $I \in \I$, there
 is a $\xi_I : \xi_{i'} \circ d(I) \Longrightarrow \xi_i$ and 
for every $\iota : I \Longrightarrow I'$, we demand $\xi_{I'} \circ_h  d\iota = \xi_I $
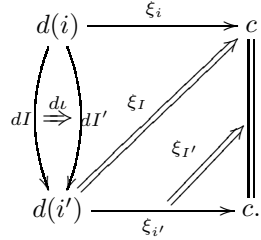
$$\xymatrix{
d(i) \ar[rr]^{\xi_i}  \ddtwocell_{dI}^{dI'}{^d\iota}
& & c \ar@{=}[dd]
\\
 & &{ }
\\
d(i') \ar@{=>}[uurr]^{\xi_I} \ar[rr]_{\xi_{i'}} & { }\ar@{=>}[ur]^{\xi_{I'}}& c.
}$$
Let $\xi'$ be another cocone over $d$ with vertex $c$, then a morphism 
of cocones $\Xi : \xi \longrightarrow \xi'$ is a family of 2-cells $\Xi_i : \xi_i \Longrightarrow \xi_i'$
indexed by the 0-cells of $\I$. These 2-cells must satisfy
$$ \Xi_i \circ_h \xi_I \quad =\quad  \xi'_{I'} \circ_v (\Xi_{i'} \circ_h d(\iota))
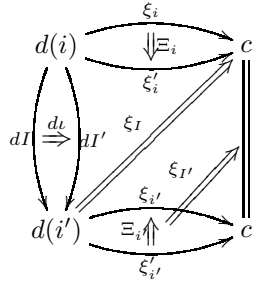
$$
$$\xymatrix{
d(i)   \ddtwocell_{dI}^{dI'}{^d\iota} \rrtwocell^{\xi_i}_{\xi'_i}{\Xi_i}
& & c \ar@{=}[dd]
\\
 & &{ }
\\
d(i') \ar@{=>}[uurr]^{\xi_I} \rrtwocell^{\xi_{i'}}_{\xi'_{i'}}{^\Xi_{i'}} & { }\ar@{=>}[ur]^{\xi_{I'}}& c
}$$
(This identity is nothing more than the cube relation  (\ref{XY:qcube})  with $F(\alpha) = d(\iota),\quad
F'(\alpha) = Id_c,\quad  \sigma_c = \xi_I, \quad \sigma'_{f'} = \xi'_{I'},\quad  \Sigma_n = \Xi_i$ and $ \Sigma_m =
\Xi_{i'}$). 

There is a need to generalize this definition. Let $\gamma:\I' \longrightarrow \I$ be a 
2-functor. $\CIg$ has the same 0-cells as $\CI$,
however, $\CIg(d,d')$ is the full subcategory of $\CI(d,d')$ consisting of 
quasi-natural transformations $\xi:d \longrightarrow d'$  such that 
$\gamma \circ \xi$ is a strict natural transformation.
There is also a generalization of $\Delta$ to $\Delta_\gamma:\C \longrightarrow \CIg$.
In detail $\CIg(d,\Delta_\gamma(c))$ consists of cocones where $\xi_{\gamma(I')}$ is
the identity for all $I' \in \I'$.

Let us move on 
to quasi-colimits \cite{Bunge1, Bunge2, Gray} of 2-diagrams. $qcolim: \CI \longrightarrow \C $ is a 2-functor
that is left $\Cat$-adjoint to $\Delta$. That is, there is an isomorphism of categories
$$\C (qcolim(d), c') \cong \CI (d,\Delta c').$$
In detail, a quasi-colimit of a diagram $d$ is a pair $(qcolim(d), \xi)$ where $qcolim(d)$ is a 0-cell
of $\C $ and $\xi$ is a cocone over $d$ with vertex $qcolim(d)$ that satisfies the following 
universal property: For any cocone over $d$ with vertex $\Delta c'$, $\xi': d \longrightarrow c'$
there is a unique $\tilde{\xi'} : qcolim(d) \longrightarrow c'$ such that $\tilde{\xi'} \circ \xi_i =
\xi_i'. $ For any $\Xi : \xi' \Longrightarrow \xi''$ there is a unique $\tilde{\Xi}:\tilde{\xi'}
 \Longrightarrow \tilde{\xi''}$ such that $\tilde{\Xi} \circ_h \xi_i = \Xi_i.$

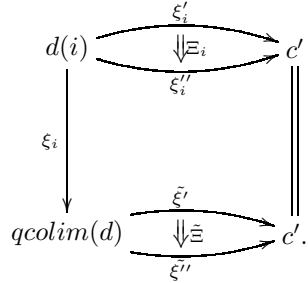
$$\xymatrix{
d(i)   \ar[dd]_{\xi_i}
\rrtwocell^{\xi'_i}_{\xi''_i}{\Xi_i} 
& & c' \ar@{=}[dd]
\\
\\
 qcolim(d) 
\rrtwocell^{\tilde{\xi'}}_{\tilde{\xi''}}
{\tilde{\Xi}} 
&& c'.
}$$

Let $\gamma:\I' \longrightarrow \I$ then $qcolim_\gamma$ is the left $\Cat$-adjoint
to $\Delta_\gamma$. In detail, we insist that $\xi_I$ in 
$$\xymatrix{
d(i)   \ar[dd]_{d(I)} \ar[rr]^{\xi_i} 
& & qcolim_\gamma(d) \ar@{=}[dd]
\\
\\
d(i') \ar[rr]_{\xi_i} \ar@{=>}[rruu]^{\xi_I} && qcolim_\gamma(d) 
}$$
be the identity if $\gamma(I') = I$ for some $I' \in \I'$.

A 2-category is {\sl quasi-cocomplete } if it has all quasi-colimits. For example,
 $\Cat$ is quasi-cocomplete. A model for $qcolim(d)$ is $\pi_0(1\downarrow d)$
where $\pi_0$ is the functor from ${\bf 2Cat}$ to $\Cat$ that forces all 2-cells to become
identities and $(1\downarrow d)$ is a 2-comma category (see page 29 of \cite{Gray} or see 
below for a definition of a 2-comma category in another context.) 
The proof that this is a model for $qcolim(d)$ is similar to the one-dimensional  colimit
case.  (Gray \cite{Gray} (page 210) proves a marvelous theorem that says that $\Cat$
is the quasi-cocompletion of $\Set$!)  If $\C$ is quasi-cocomplete, then for all 
$\T 2$, and for all $G_2:\T 1 \longrightarrow \T 2$,   $\TAC 2$ is
 also quasi-cocomplete since one can put a $\T 2$ structure on the quasi-colimit.

A {\sl weak-terminal object} of a 2-category $\I$ is a 0-cell $t\in \I$
with the following property: for every 0-cell $i \in \I$ there 
is a 1-cell $l:i \lra t$ and for any two 1-cells $l,l':i \lra t$
there is a unique iso-2-cell $\iota:l \Longleftrightarrow l'$.
Let $\gamma : \I' \lra \I$ be a 2-functor. A $\gamma$-relative
terminal object is a weak terminal object in $\I$ with the 
added requirement that if $i=\gamma(i')$ for some $i' \in \I'$
then $\iota$ is the identity i.e. $l$ is a {\sl unique} 1-cell.
If $\gamma=id_{\I}$ then a $\gamma$-relative terminal object
is, in fact, a terminal object. If $\gamma$ is the unique 2-functor from the 
empty 2-category to $\I$ then a $\gamma$-relative terminal object
is a weak-terminal object.  
Whereas a terminal object is unique up to a  
unique isomorphism, a weak-terminal object is unique up to an equivalence.
To see this, let $t_1$ and $t_2$ be weak-terminal objects. We then have
$$\xymatrix{ t_1 \ar[r]^f \ar@/_2pc/[rr]_{id}^{\Updownarrow}
 & t_2 \ar[r]^g\ar@/^2pc/[rr]^{id}_{\Updownarrow} & t_1 \ar[r]^f & t_2. }$$
The reader inclined to think topologically should think of the 
terminal object as the one-point topological space and a weak-terminal 
object as a contractible pointed space.

\begin{prop} 
Let $t\in \I$ be a $\gamma$-relative terminal object and let $d:\I\lra\C$ be a 2-diagram, then $qcolim_\gamma(d)$
is equivalent to $d(t)$ . If $\gamma = id_\I$ then $qcolim_\gamma(d)$
is isomorphic to $d(t)$ 
\end{prop}

Given $G:\TO 1 \longrightarrow \TO 2$ and a 0-cell $n \in \T 2$,
we define  the 2-comma category $(G\downarrow n)$.  0-cells are 
pairs $(Gm, g:Gm \longrightarrow n)$; 1-cells are pairs $(Gh:Gm \longrightarrow
Gm', \tau_h: g' \circ Gh \Longrightarrow g)$  where $\tau_h$ is a 2-cell in $\TO 2$
that makes 
$$ \xymatrix{
G(m) \ar[dd]_{Gh} \ar[rr]^g && n\ar@{=}[dd]\\
 & &{ }\\
G(m') \ar[rr]_{g'} \ar@{=>}[uurr]^{\tau_h}&&n
}$$
commute; 2-cells are $G(\beta):Gh \Longrightarrow Gh'$ that satisfy
$\tau_h = \tau_{h'} \circ G(\beta)$.  If $G_2: (\TO 0)' \lra \TO 2$ then 
$(G \downarrow n)_{G_2}$ is the locally full subcategory where $\tau_h = id$
if $h=G_2(h')$ for some $h'$ in $(\TO 0)'$.
$f:n \longrightarrow n'$ in $\TO 2$
induces a 2-functor $(G\downarrow f):(G\downarrow n) \longrightarrow (G\downarrow n')$.
$\alpha : f \Longrightarrow f'$ in $\TO 2$ induces a 2-natural transformation 
$(G\downarrow \alpha) : (G\downarrow f) \Longrightarrow (G\downarrow f')$. There is also an obvious forgetful
2-functor $P:(G\downarrow n) \longrightarrow \TO 1$ that commutes with $(G\downarrow f)$ 
and $(G\downarrow \alpha)$. There are similar properties for $(G \downarrow n)_{G_2}$.

The final preliminary needed is 
\begin{defi}
Let $\A$ and $\B$ be 2-categories. Let $L:\A \lra \B$ and 
$G:\B \lra \A$ be 2-functors. $L$ is a (strict c.f. pg 168 
of \cite{Gray}) left quasi-adjoint of $G$ if there exists two 
quasi-natural transformations $\eta:id_\A \lra GL$ and 
$\varepsilon:FG \lra id_\B$ strictly satisfying
the usual two triangle identities.
\end{defi}

Every $G:\TO 1 \longrightarrow \TO 2$ within  a commutative square 
\be \xymatrix{\label{XY:2thsquare}
\TO 0 \ar[dd]_{G_0} \ar[rr]^{G_1} && \TO 1 \ar[dd]^{G}
\\
\\
(\T 0')^{op} \ar[rr]_{G_2} && \TO 2
}\ee
induces a 2-functor $G^*: {\bf Hom}(\TO 2, \C) \longrightarrow  {\bf Hom}(\TO 1, \C)$
via precomposition. From the fact that $G$ preserves products and the square (\ref{XY:2thsquare})
commutes, $G^*$
restricts to an algebraic 2-functor 
$G^*:\TAC 2 \longrightarrow \TAC 1.$

\begin{teo}
Let $\C$ be a Cartesian closed quasi-cocomplete 2-category.
 Every $G^*:\TAC 2 \longrightarrow \TAC 1.$ has a
(strict) left quasi-adjoint $Lan_G(F):\TAC 1 \longrightarrow \TAC 2$
which can be computed ``pointwise'' for $F$ in $\TAC 1$ as 
$$Lan_G(F)(?)  =  qcolim_{G_1}(F\circ P: (G \downarrow (?))_{G_2} \longrightarrow \TO 1 \longrightarrow
\C).$$ Furthermore: i) $Lan_G$ takes quasi-natural transformations to natural
transformations. ii)$\eta_F$ has a left inverse and iii) $\varepsilon_K$ has a 
right inverse.
\end{teo}
{\bf Proof.}
By cocompleteness of $Lan_G(F)(n)$ is an object of $\C$.
We prove the many steps in small bites:

\underline{ $Lan_G(F)$ is a 2-functor.} $f:n \longrightarrow n'$
induces  $(G\downarrow f)_{G_2}:(G\downarrow n)_{G_2} \longrightarrow (G\downarrow n')_{G_2}$.
which induces $\tilde{f}:Lan(F)(n) \longrightarrow Lan(F)(n')$.
We must stress that for all $g:Gm \longrightarrow n$ in $(G\downarrow n)_{G_2}$
\be \xymatrix{\label{XY:functofLan}
&& Lan(F)(n) \ar[dd]^{\tilde{f}}
\\
Fm\ar[urr]^{\xi_g} 
\ar[drr]_{\xi_{f\circ g}}
\\
&& Lan(F)(n')}
\ee
commutes strictly. There is a similar picture for $\alpha :  f \Longrightarrow f'$
and the induced $\tilde{\alpha}:Lan(F)(f) \Longrightarrow Lan(F)(f')$ .

\underline{$Lan_G(F)$ preserves products.} This is very similar to the one dimensional
case and we leave it for the reader. Preservation of products is not true for 
all 2-categories but it is true for 
the usual 2-categories that one takes algebras in, like $\Cat$ and any other Cartesian 
closed 2-category.

\underline{$Lan_G$ takes quasi-natural transformations to natural transformations.}
Let $\sigma: F \longrightarrow F'$ be a quasi-natural transformation in $\TAC 1$.
$\sigma_m:Fm \longrightarrow F'm$  makes $Lan(F')(n)$ satisfy the universal property 
of $Lan(F)(n)$  and  so we have the commuting square
\be \xymatrix{\label{XY:functofLan2}
Fm \ar[rr]^{\sigma_m} \ar[dd]_{\xi}
&&
F'm \ar[dd]^{\xi'}
\\
\\
Lan(F)(n) \ar@{-->}[rr]^{\tilde{\sigma_n}}
&& Lan(F')(n).}\ee

For $f:n\longrightarrow n'$, there is 
$$\xymatrix{
Fm \ar[rrrr]^{\sigma_m} \ar@{=}[dddd] \ar[dr]
&&&&
F'm \ar@{=}[dddd]\ar[dl]
\\
& Lan(F)(n) \ar@{-->}[rr]^{\tilde{\sigma_n}} \ar[dd]_{LanFf}
&& Lan(F')(n)\ar[dd]^{LanF'f}
\\
\\
& Lan(F)(n') \ar@{-->}[rr]^{\tilde{\sigma_{n'}}} 
&& Lan(F')(n')
\\
Fm \ar[rrrr]^{\sigma_m} \ar[ur]
&&&&
F'm \ar[ul]}$$ 
where the left and right quadrilaterals commute from diagram (\ref{XY:functofLan}).
 The top and bottom quadrilaterals commute from diagram (\ref{XY:functofLan2}).
 Since $Lan(F')(n')$ satisfy the universal properties of $Lan(F)(n)$ there is a 
unique  $Lan(F)(n) \longrightarrow Lan(F')(n')$ which coheres with the 
surrounding commutative quadrilaterals. Hence the inner square commutes
making $\tilde{\sigma}$ a natural transformation ({\em not quasi}) in $\TAC 2$.

\underline{The unit of the $G^*\vdash Lan_G$ quasi-adjunction:
 $\eta_F:F \longrightarrow (G^*\circ Lan_G)(F)$}
$$(G^*\circ Lan_G)(F) = Lan_G(F)\circ G =  qcolim_{G_1}(F\circ P: (G \downarrow G(?)_{G_2}) 
\longrightarrow \TO 1 \longrightarrow \C).$$
Within $(G\downarrow G(n))_{G_2}$ there is $id:G(n) \longrightarrow G(n)$.
$(F\circ P)(id:G(n) \longrightarrow G(n))=Fn$. And so we set 
$$[\eta_{F,n} = \xi_{id}]:Fn \longrightarrow [qcolim_{G_1}= (G^*\circ Lan_G)(F)(n)].$$
Given $f:n \longrightarrow n'$ in $\TO 1$, we have $G(f):Gn \longrightarrow Gn'$ in 
$\TO 2$ which induces $(G\downarrow G(f))_{G_2}:(G\downarrow Gn)_{G_2} \lra (G\downarrow Gn')_{G_2}$ and 
hence $\tilde{Gf}:G^*Lan(F)(n) \longrightarrow G^*Lan(F)(n')$.
Which makes $\eta_F$ a quasi-natural transformation:
$$ \xymatrix{
Fn \ar[rr]^{\eta_{F,n}} \ar[dd]_{Ff}\ar[rrdd]^{\xi_{G(f)}}&& G^*Lan(F)(n) \ar[dd]^{G^*Lan(F)(f)}\\
&{ }&{ }\\
Fn' \ar[rr]_{\eta_{F,n'}} \ar@{=>}[ru]^{\xi_f}& & G^*Lan(F)(n') \\}
$$ 
where the upper right triangle commutes.

\begin{remar} \label{natunit}
If $f=G_1(\bar{f})$ for some $\bar{f}$ in $\TO 0$ then $\xi_f$ is the 
identity making the square commute. If this is true for all $f \in \TO 1$ then $\eta_F$ is, in fact,
a natural transformation.
\end{remar}

One gets the left inverse of $\eta_F$ from the commutativity of the 
bottom quadrilateral of
$$\xymatrix{
F(m) \ar[rr]^{id} \ar[dd] \ar[dr]&& F(n)n \ar@{=}[d]
\\
&G^*LanFn \ar@{-->}[r]^{\beta}&F(n)\ar@{=}[d]
\\
F(n) \ar[ur]^{\eta_{F,n}} \ar[rr]_{id} &{ }\ar@{=>}[ruu]|\hole & F(n).}
$$

\underline{The counit of the $G^*\vdash Lan_G$ quasi-adjunction:
 $\varepsilon_K: (Lan_G \circ G^*)(K) \lra K$}
$$ (Lan_G \circ G^*)(K) = Lan_G(KG) = qcolim_{G_1}(K\circ G \circ P: (G \downarrow (?))_{G_2} 
\longrightarrow \TO 1 \longrightarrow \TO 2 \longrightarrow \C).$$
Consider the following typical diagram in $(G\downarrow n)$
$$\xymatrix{
G(m) \ar[rr]^g \ar[dd]_{Gh}&& n\ar@{=}[dd]
\\ &&{ }
\\
G(m') \ar[rr]_{g'} \ar@{=>}[rruu]^{\tau_h}&& n.}
$$
Applying $qcolim_{G_1}(K\circ G \circ P)$ and $K$ to this diagram
gives us:
$$\xymatrix{ \label{XY:counitsquare}
KG(m) \ar[rr]^{Kg} \ar[dd]_{KGh}\ar[dr]&& Kn \ar@{=}[d]
\\
&LanG^*Kn \ar@{-->}[r]^{\varepsilon}&Kn\ar@{=}[d]
\\
KG(m') \ar[ur] \ar[rr]_{Kg'} &{ }\ar@{=>}[ruu]|\hole_{K\tau_h}& Kn.}
$$
And so there is the induced $\varepsilon_{K,n}: (Lan_G \circ G^*)(K)(n) \lra Kn$.

$\varepsilon_K$ is also quasi-natural. Given $f:n \lra n'$
in $\TO 2$ we have the diagram in $(G \downarrow n')$
$$\xymatrix{
G(n) \ar[rr]^{id} \ar[dd]_{G(\hat{f})}&& n \ar[dd]^f
\\
&& { }
\\
G(n') \ar[rr]_{id} \ar@{=>}[rruu]&& n'
}
$$
where $\hat{f}$ is in $\TO 1$. This square commutes
if $\hat{f} = G_1(\bar{f})$ for some $\bar{f} \in \TO 0$.
Applying $K$ and taking appropriate quasi-colimits
to this commutative or noncommutative diagram
gives us:
$$\xymatrix{
KG(n) \ar[rrr]^{id} \ar[dr] \ar[dddd]_{KG(\hat{f})} &&&K(n) \ar@{=}[d]
\\
& Lan G^*Kn \ar[rr]^{\varepsilon_{K,n}} \ar[dd]_{LanGKf}
\ar[rrdd] &&K(n) \ar[dd]^{K(f)}
\\
&& { }\ar@{=>}[ru]
\\
& Lan G^*Kn' \ar[rr]^{\varepsilon_{K,n'}} && K(n') \ar@{=}[d]
\\
KG(n') \ar[ru] \ar[rrr]_{id} &&&K(n')
}
$$
where the surrounding quadrilateral and the lower left triangle
commute. If the outer square commutes then the inner square must also
commute.

\begin{remar} \label{natcounit}
If $\hat{f} = G_1(\bar{f})$ for some $\bar{f} \in \TO 0$, then the
outer square commutes. If this is true for all $f \in \TO 2$
then $\varepsilon_K$ is, in fact, a natural transformation.
\end{remar}

One gets the right inverse of $\varepsilon_K$ from the commutativity
of the bottom quadrilateral of the diagram that gives $\varepsilon_K$ 
when you
set $m' = n$ and $g'= id_n$.

We leave the following usual two triangle identities for the reader's pleasure:
$$\xymatrix{
G^*Kn \ar[dd]_{\eta G^*K} \ar[rrdd]^{id}
&&
Lan(F)n \ar[rr]^{Lan \eta_F} \ar[rrdd]_{id}
&&
Lan \circ G^* \circ Lan(F)n \ar[dd]^{\varepsilon_{LanF}}
\\
\\
G^* \circ Lan \circ G^* Kn \ar[rr]_{G^* \varepsilon}
&&
G^*Kn
&&
Lan(F)n}
$$
Q.E.D.

Now that we have these tools, we can go on and prove coherence theorems.
Following Remark.3.1 (respectively Remark.3.2 ) we have

\begin{teo}
If there exists a 2-theory-morphism $H:\TO 2 \lra \TO 0$ (resp. 
$H':\TO 1 \lra (\TO 0)'$ such that 
$$\xymatrix{
\TO 0 \ar[rr]^{G_1} \ar[dd] && \TO 1 \ar[dd]^{G} 
&&
\TO 0 \ar[rr]^{G_1} \ar[dd] && \TO 1 \ar[dd]^{G} \ar[ddll]_{H'} 
\\
&& & \mbox{(resp.} &&&&\mbox{)}
\\
(\TO 0)' \ar[rr]^{G_2}  && \TO 2 \ar[lluu]_{H}
&&
(\TO 0)' \ar[rr]^{G_2}  && \TO 2  
}
$$
the two triangles commute, then the unit (resp. counit) of 
the $G^* \vdash Lan_G$ adjunction is a natural transformation.$\Box$
\end{teo}

Setting $\TO 0 = \TO 1 =\finbar^{op}$ gives us the naturality of the 
counit and so we have

\begin{cor}
For  any $G_2:\TO 0 \lra \TO 2$ we have the following isomorphism of
2-categories:
$${\bf 2Alg_{G_2}(T_2, Cat)}(Lan_G F, K) \quad \cong \quad \Cat (F, G^*K)$$
where $F$ is a category (functor from the trivial theory) and $K$ is a $\TT_2$-algebra. i.e. $Lan_G F$ is the free $\TT_2$ category over F.
\end{cor}
{\bf Proof.} The only non-obvious part is the universality of the 
counit. This is similar to the one dimensional case. We leave the following
diagram to help:

$$\xymatrix{
F(n) \ar[rr]^{\xi} \ar@{-->}[rrdd] \ar@{-->}[dd]&&
LanF(n) \ar@{-->}[dd] \ar[rr]^{\alpha} &&Kn \ar@{=}[dd]
\\
\\
KG(n)\ar[rr]^{\xi} && LanKG(n) \ar[rr]_{\varepsilon_{K,n}} && Kn.\quad \Box
}
$$

Setting $\TO 0$ to also be $\finbar^{op}$ gives us an unrestricted
$Lan F$. This is used in the reconstruction of a theory from its 
category of algebras (For technical reasons from the quasi-Yoneda lemma 
we insist that the category of algebras have quasi-natural transformations
where the squares commute up to a {\bf iso}-2-cell.)

\begin{teo}
Every theory $\TT$ is quasi-equivalent to its 2-category of algebras, ${\bf 2Alg^i(T, Cat)}$.
\end{teo}
{\bf Proof} Let $F_{[n]}:\finbar^{op} \lra \Cat$ be the ``constant''
functor on $[n]$ i.e. $F_{[n]}(m) \cong [n]^m$. $Lan F_{[n]}$ is the 
free $\TT$-algebra on $[n]$ elements. By Proposition 2, $LanF_{[n]} \cong
\TT^{op}(n,-)$. Using this and the quasi-Yoneda lemma (e), we get the
following quasi-equivalence of categories 
$$\begin{array}{rcl}
{\bf 2Alg^i}(\TT, \Cat)( F_{[n]},F_{[m]}) &\cong&{\bf qNat^i}(F_{[n]},F_{[m]})\\
&\cong& {\bf qNat^i}(\TT^{op}(n,-),\TT^{op}(m,-))\\
&\cong& \TT^{op}(n,m)\qquad \Box.
\end{array}$$  
This quasi-equivalence is the counit of the adjunction:
$$\xymatrix{
\TTh  
\ar@<2ex>[rrrr]^{2Alg(-,Cat)}_{\bot}
&&&&\TCatcat .
\ar@<2ex>[llll]^{\mbox{Free}}
}$$
 Where $\TCatcat$ denotes the tractable 2-functors $U:\C \lra
\Cat$. Tractable means $\C$ must be a local groupoid and the category 
${\bf 2Cat}(U^n, U^m)$
must be small and locally small. There will be more about this adjunction
at the end of section 4.

\begin{defi} A weakly-unique quasi-section of $G:\TO 1 \lra \TO 2$
is a 2-theory-morphism $H:\TO 2 \lra \TO 1$ satisfying:
\begin{enumerate}
\item the diagram 
$$\xymatrix{
\TO 0 \ar[rr]^{G_1} \ar[dd] && \TO 1
\\
&& 
\\
(\TO 0)' \ar[rr]_{G_2}  && \TO 2 \ar[uu]_H
}
$$
commutes
\item for every 1-cell $f \in \TO 2$ there is a 2-cell $\alpha:(G \circ H)(f)
\Longrightarrow f$
\item $H$ is unique up to a {\bf unique} 2-cell.
\end{enumerate}
\end{defi}

\begin{teo}
If $G$ has a weakly-unique quasi-section, then $\eta_F$ is an equivalence
for every $F \in {\bf 2Alg_{G_1}}(\TO 1, \Cat)$.
\end{teo}
{\bf Proof.}
The definition of a weakly-unique quasi-section insures that $id:G(n) \lra G(n)$ is a $G_1$-relative terminal object of $(G\downarrow G(n))$. Hence, 
using Proposition 3, we have that $\eta_F$ is an equivalence. $\Box$

\begin{examp}
Setting $\TO 0 = (\TO 0)' =\finbar^{op}$ and $G:\Tmon \lra \Tsmon$ (the obvious
``strictification'' functor), we have that every monoidal category is tensor
equivalent to a strict monoidal category.
\end{examp}

\begin{examp}
Following the above, we have a more general theorem. 
Let $\TO X$ be any theory that ``contains''
the monoidal theory. 
Let ${\bf T^{op}_{sX}}$ be the strict version of that theory with 
$G$ being the ``strictification'' 2-theory-morphism. 
Then $\eta_F$ is an equivalence.
\end{examp}

As with all conditions, the case where a condition fails is far 
more interesting. For example $G:\Tassoc \lra \Tsmon$ has many quasi-sections 
but they are 
not unique up to a {\bf unique} isomorphism. Similarly for $G:\Tbraid 
\lra \Tsym$.  Notice that in all these cases, $Lan_GF$ always 
exist and there are many things that one can say about $\eta_F$. But it is
not an equivalence. There is much structure to explore.

Many other coherence theorems can be stated and proved on the syntactical 
level. For example, Corollary 2.4 (pg. 43) of \cite{Joy&Str} says that 
given 
$$\xymatrix{
\finbar^{op} \ar[drr]^G \ar[dd]
\\
&& \Tsbraid^{op}
\\
\Tbraid^{op} \ar[rru]_{G_1}}
$$
$\eta_F:F \lra G^* Lan_GF$ determines an equivalence of categories
$${\bf 2Alg}(\Tbraid, \Cat)(G_1^*(Lan_GF), G_1^*(V)) \cong \Cat(F, G^*(V))$$
This is proven using the properties of $G, G_1$ and the universal properties
of $Lan_G$.

One can go on to formalize many coherence statements like ``If $G$ is 
locally faithful etc ... and $G_1$ is faithful etc ... , then
$\eta_F$ is ... and $\varepsilon_K$ is ....'' . We leave this 
noble task for future explorers.

\section{Kronecker product}
It is common to look at the algebras of one theory in the category of algebras 
of another theory. The theory of such algebras is given as the Kronecker
product of the two theories. 

The Kronecker product \cite{Freyd} of (1-)theories is a well understood
coherent symmetric monoidal 
2-bifunctor $\otimes_K:\th \times \th \longrightarrow \th$. 
Let $T_1$ and $T_2$ be two theories. $T_1 \otimes_K T_2$ is a theory that satisfies
the  universal property
$$ Alg(T_1 \otimes_K T_2, C) \cong  Alg(T_1, Alg(T_2,C)).$$
$T_1 \otimes_K T_2$ is constructed as follows. Construct the 
the coproduct in the category of theories (pushout in $\Cat$)
$$\xymatrix{
\fin \ar[rr] \ar[dd]&&T_1 \ar[dd]
\\
\\
T_2 \ar[rr]&& T_1 \coprod T_2.
}$$
Place a congruence on  $ T_1 \coprod T_2$ such  that for all 
 $f:m\longrightarrow m'$ in $T_1$ and $g:n \longrightarrow n'$  in 
$T_2$
the diagram 
$$ \xymatrix{
		&n^m \ar[rr]^{g^m} &&n'^m \ar[dr]^{\sim}
\\
m^n \ar[dd]_{f^n}\ar[ur]^{\sim} &&&&m^{n'} \ar[dd]^{f^{n'}}
\\
\\
{m'}^n \ar[dr]^{\sim}&&&&{m'}^{n'} 
\\
		&n^{m'} \ar[rr]^{g^{m'}} &&{n'}^{m'} \ar[ur]^{\sim}
}$$
commutes. We have a  full theory-morphism
$ T_1 \coprod T_2 \longrightarrow T_1 \otimes_K T_2.$

There is an analogous Kronecker product on the semantic level. Although
we have not been able to find this construction in the literature, surely it is
well known to the cognoscenti. Denote the tractable 2-functors from
 $\Cat$ to $\Set$  as $\catset$. The semantic Kronecker product 
is a coherently monoidal symmetric 2-bifunctor
$\oplus_K:\catset \times \catset \longrightarrow \catset$. Given two such
tractable functors $U_1: C_1 \longrightarrow \Set$
and 
$U_2:C_2 \longrightarrow \Set$, $C_1 \oplus_K C_2 \longrightarrow \Set$
is constructed as follows. Construct the product in $\catset$
 (pullback in $\Cat$)
 $$\xymatrix{
C_1 \times_{\Set}C_2 \ar[rr] \ar[dd]&&C_1 \ar[dd]
\\
\\
C_2 \ar[rr]&& \Set .
}$$
$C_1 \times_{\Set}C_2$ is to be thought of as sets with
both a  $C_1$ structure and a $C_2$ structure.  $C_1 \oplus_K C_2$
is the full subcategory of $C_1 \times_{\Set}C_2$ consisting of those 
objects $c$ that satisfy the following condition: for all $f:c^m \longrightarrow
 c^{m'} $in $C_1$ and $g:c^n \longrightarrow c^{n'}$ in $C_2$,
$$ \xymatrix{
		&U_2(c^n)^m \ar[rr]^{U_2(g)^m} &&U_2(c^{n'})^m \ar[dr]^{\sim}
\\
U_1(c^m)^n \ar[dd]_{U_1(f)^n}\ar[ur]^{\sim} &&&&U_1(c^m)^{n'} \ar[dd]^{U_1(f)^{n'}}
\\
\\
U_1(c^{m'})^n \ar[dr]^{\sim}&&&&U_1(c^{m'})^{n'} 
\\
		&U_2(c^n)^{m'} \ar[rr]^{U_2(g)^{m'}} &&U_2(c^{n'})^{m'} \ar[ur]^{\sim}
}$$
commutes. It is not hard to show that the structure - semantics adjunction 
(equivalence) takes the Kronecker product theories to the 
Kronecker product semantics and vice versa. See the end of section 4 for a large diagram 
showing what commutes.

There is a two-dimensional analogue to the Kronecker product. Rather than look
at two 2-theories $\T 1$ and $\T 2$ that are disconnected, we shall assume that 
both of these theories have an underlying $\T 0$, i.e. there is a diagram in $\TTh$
$$\xymatrix{ 
\T 1 &&\T 0 \ar[ll]_{G_1} \ar[rr]^{G_2}&&\T 2.}$$
We can, however, give a similar  bifunctor
 which assumes no underlying  $\T 0$ (i.e. $\T 0 = \finbar$)
 or assuming a commutative 
 square of 2-theories: 
$$\xymatrix{
 \T 0 \ar[rr]^{G_1}\ar[dd]&&\T 1 \ar[dd]
 \\
 \\
 \T 0 ' \ar[rr]^{G_2}&& \T 2.
}$$
However most examples are found with one  underlying 2-theory.

\begin{defi} 
A (2-)Kronecker product of 2-theories is a 3-bifunctor
$$\KP :(\T 0 \downarrow \TTh )\times (\T 0 \downarrow \TTh) \longrightarrow
(\T 0 \downarrow \TTh)$$
that satisfies the following universal property: for all 
$$\xymatrix{\T 1 & & \T 0 \ar[ll]_{G_1} \ar[rr]^{G_2} && \T 2}$$
there is an induced 
$$\xymatrix{
\T 0 \ar[rr]^{G_1}\ar[dd]_{G_2}\ar[ddrr]^{G_1 \KP G_2} &&\T 1 \ar[dd]
\\
\\
\T 2 \ar[rr]&& \T 1 \KP \T 2.
}$$
and for all 2-categories with finite products $\C$, an isomorphism  of 2-categories
\be {\bf 2Alg}_{G_1\KP G_2}(\T 1 \KP \T 2, \C) \cong {\bf 2Alg}_{G_1}(\T 1, {\bf 2Alg}_{G_2}(\T 2, \C))\label{KPuniv}\ee
which is natural for all cells in $(\T 0 \downarrow \TTh)$ and for all cells in $\C$.
\end{defi}

When $\C$ is ``nice'' and the 2-theory is reconstructible from 
its 2-category of algebras we have  
$$\begin{array}{rcl}
{\bf 2Alg}(\T 1 \KP (\T 2 \KP \T 3), \C) &\cong&{\bf 2Alg}(\T 1, {\bf 2Alg}(\T 2 \KP \T 3, \C)) \\
								 &\cong&{\bf 2Alg}(\T 1, {\bf 2Alg}(\T 2, {\bf 2Alg}(  \T 3, \C))) \\
&\cong&{\bf 2Alg}((\T 1\KP \T 2), {\bf 2Alg}(  \T 3, \C)) \\
&\cong&{\bf 2Alg}((\T 1\KP \T 2) \KP \T 3, \C). 
\end{array}$$  
and hence $\T 1 \KP (\T 2 \KP \T 3) \cong (\T 1\KP \T 2) \KP \T 3.$ It is conjectured that
this bifunctor is actually coherently associative (c.f. \cite{Grayc}) but we leave this 
question for now. If we insist that the Kronecker product satisfy
$${\bf 2Alg}_{G_1\KP G_2}^{\bf i}(\T 1 \KP \T 2, \C) \cong {\bf 2Alg}_{G_1}^{\bf i}(\T 1, 
{\bf 2Alg}_{G_2}^{\bf i}(\T 2, \C))$$
then $\T 1 \KP \T 2$ will be (coherently) isomorphic to $\T 2 \KP \T 1$.

In order to construct $\T 1 \KP \T 2$, we take the coproduct 
 $\T 1 \coprod_{\T 0} \T 2$ in $(\T 0 \downarrow \TTh)$ and 
we freely add in the following 2-cells: For every $f:m \longrightarrow m'$
in $\T 1$ and $g:n \longrightarrow n'$ in $\T 2$  we add the 2-cell $\delta 
(f,g)$ that makes the following diagram commute:
\be  \xymatrix{\label{XY:KPquasi}
		&n^m \ar[rr]^{g^m} &&n'^m \ar[dr]^{\sim}
\\
m^n \ar[dd]_{f^n}\ar[ur]^{\sim} &&&&m^{n'} \ar[dd]^{f^{n'}}
\\
&&&&{ }
\\
{m'}^n \ar[dr]^{\sim}&&&&{m'}^{n'} 
\\
		&n^{m'} \ar[rr]^{g^{m'}} \ar@{=>}[rrruuu]^{\delta 
(f,g)}&{ }&{n'}^{m'} \ar[ur]^{\sim}
}\ee
[If $\KP$ is to be symmetric, then we must insist that $\delta (f,g)$
be an isomorphism.]
The $\delta$`s must satisfy the following coherence conditions that are 
compatible to the four coherence conditions in the definition of
a quasi-natural transformation.
\begin{enumerate}
\item If $f$ is in the image of $G_1$ [or if $g$ is in the image of $G_2$], then
$\delta (f,g)$ must be set to the identity.
\item $\delta$ must preserve products in $f$ [and $g$] as in Figure I.
\item  $\delta (f\circ f',g) = \delta(f,g) \circ_v \delta(f',g) $ [and 
  $\delta (f, g\circ g') = \delta(f,g) \circ_h \delta(f,g') $.]
 \item $\delta$ must preserve 2-cells. i.e. If there is a 2-cell in $\T 1$
$$ 
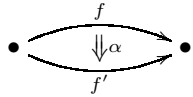
\xymatrix{\bullet \rrtwocell^f_{f'}{\alpha} && \bullet
}$$
then we have the following equality of diagrams (we leave out the corner isomorphisms and the exponents) 
$$\xymatrix{
\bullet \ar[rr]^g \ar[dd] \ddlowertwocell_f^{f'}{^\alpha}& & \bullet \ar[dd]^{f'} & &
\bullet \ar[rr]^g \ar[dd]_{f} & & \bullet  \ar[dd]\dduppertwocell^{f'}_f{^\alpha}
\\
 & &{ } &=&&& { }
\\
 \bullet \ar[rr]_g &{ }\ar@{=>}[ru]^{\delta(f',g)} & \bullet & &
\bullet \ar[rr]_g &{ }\ar@{=>}[ru]^{\delta(f,g)}  & \bullet.
 }$$
[For symmetry, a 
$$ \xymatrix{\bullet \rrtwocell^g_{g'}{\beta} && \bullet
}$$
in $\T 2$
implies  
$$\xymatrix{
\bullet\rruppertwocell^{g'}_g{^\beta} \ar[rr] \ar[dd]_{f} & & \bullet \ar[dd]^{f} & &
\bullet \ar[rr]^{g'} \ar[dd]_{f} & & \bullet \ar[dd]^{f}
\\
 & &{ } &=&&& { }
\\
 \bullet \ar[rr]_g &{ }\ar@{=>}[ru]^{\delta(f,g)} & \bullet & &
\bullet \ar[rr]\rrlowertwocell_g^{g'}{^\beta} &{ }\ar@{=>}[ru]^{\delta(f,g')}  & \bullet
 .]}$$
\end{enumerate}

\begin{remar} We demand not only that  the $\delta$'s preserve 
the Cartesian product, but that the  $\delta$'s
inherit {\it all} the coherence properties
of the Cartesian product.
\end{remar}

The fact that there is choice in the construction of $\T 1 \KP \T 2$, should 
not disturb the reader too much since we never claimed that $\T 1 \KP \T 2$
should be unique. Rather, it should be unique up to a (2-)isomorphism.
In order to show that our construction of $\T 1 \KP \T 2$
satisfies the universal properties 
stated in \ref{KPuniv}, let us examine an algebra in
${\bf 2Alg}_{G_1}(\T 1, {\bf 2Alg}_{G_2}(\T 2, \C))$.
An algebra is a finite product preserving functor $F:\T 1 \lra \C$.
Assume $F(1) = G:\T 2 \lra \C$. Then 
$F(m)=F(1)^m=G^m: \T 2 \lra \C$  For every $f:m \lra m'$ in $\T 1$,
$F(f)$ is a quasi-natural 
transformation from $G^m$ to $G^{m'}$. In order for $F(f)$ to be 
such a quasi-natural transformation,
we must have that every $g:n \lra n'$ in $\T 2$, makes a the following diagram:
\be  \xymatrix{
		&G^{n^m} \ar[rr]^{g^m} &&G^{n'^m} \ar[dr]^{\sim}
\\
G^{m^n} \ar[dd]_{f^n}\ar[ur]^{\sim} &&&&G^{m^{n'}} \ar[dd]^{f^{n'}}
\\
&&&&{ }
\\
G^{{m'}^n} \ar[dr]^{\sim}&&&&G^{{m'}^{n'}} 
\\
		&G^{n^{m'}} \ar[rr]^{g^{m'}} 
\ar@{=>}[rrruuu]^{F(f)_g}&{ }&G^{{n'}^{m'}} \ar[ur]^{\sim}
}.\ee This is what is described by \ref{XY:KPquasi}.
$F(f)_g $ is what corresponds to  $\delta(f,g)$ in our theory 
$\T 1 \KP \T 2$. The rest of the tedious details melt away when
one realizes that our construction was made to mimic the definition of a
quasi-natural transformation in our 2-categories of
algebras.

There is a similar construction for the Kronecker product on the 
semantic level: 
$$
\oplus_K: \TCatcat \times \TCatcat \longrightarrow \TCatcat.$$
We leave the details for the reader.

Our examples have all been proven in Joyal and Street's paper
\cite{Joy&Str}. We are simply restating them in the language of 
Kronecker products.

All our examples have the $\delta$'s as isomorphisms. We, however,
must stress that this is an historical accident rather then something 
intrinsicly important to coherence. Even though most coherence results 
are about natural
{\em iso}morphisms, one should study the 
general case where the  
natural transformations are in not necessarily isomorphisms. The only 
example in the literature that we know of where coherence questions arise 
for natural transformations that are not isomorphisms is Yetter's 
notion of a pre-braiding $\cite{Yetter}$.

In order to make the diagrams in the examples a little more readable,
 we shall write our morphisms of the  theories the opposite way. In other
 words, we shall write them as if $\TT$ was 
$\TT^{op}$.

\begin{examp}
Let $\Tsmon$ be the theory of strict (associativity) monoidal categories. 
Let $\T 0$ be the theory of pointed categories, that is the theory of 
categories with a distinguished element to be thought of as a unit 
of the tensor product(s). Then we have 
$$ \Tsmon \KP \Tsmon \cong \Tsbraid$$
where $\Tsbraid$ is the theory of braided monoidal categories. This result
is a two dimensional version of the fact that the Kronecker product of the 
theory of monoids with itself is the theory of commutative monoids.

In order to distinguish the two (isomorphic) multiplications, we
shall denote one by $\otimes:2 \longrightarrow 1$ and one 
by $\Phi:2 \longrightarrow 1$. By the construction of the Kronecker product we have (abandoning corner isomorphisms)
$$
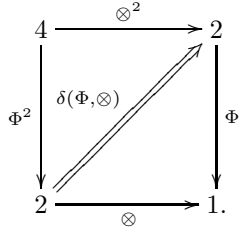
\xymatrix{
4 \ar[rr]^{\otimes^2} \ar[dd]_{\Phi^2} &&
2 \ar[dd]^{\Phi}
\\
\\
2 \ar[rr]_{\otimes}\ar@{=>}[rruu]^{\delta(\Phi,\otimes)}&&1.
}$$

On the semantic level, $\delta(\Phi,\otimes)$ induces an
isomorphism 
$$ \delta(\Phi,\otimes)_{A,A',B,B'} : (A\Phi B) \otimes (A' \Phi B')
\longrightarrow (A\otimes  A') \Phi (B \otimes B').$$
Setting $A'=B=I$ , we get an isomorphism $A\otimes B' \longrightarrow A \Phi
B'$. 
Setting $A=B'=I$ we get an isomorphism $B \otimes A' \longrightarrow A' \Phi B$. And so setting 
$$\gamma_{A,B}=\delta^{-1}_{I,A,B,I}\circ \delta_{A,I,I,B} :A\otimes B
 \longrightarrow A \Phi B \longrightarrow
B \otimes A.$$ Only the braiding relation is left to be shown. By creatively pasting the coherence 
conditions for $\delta$ (i.e. $\delta \circ (\delta \times 1) = \delta \circ (1 \times \delta)$),  we have 
$$
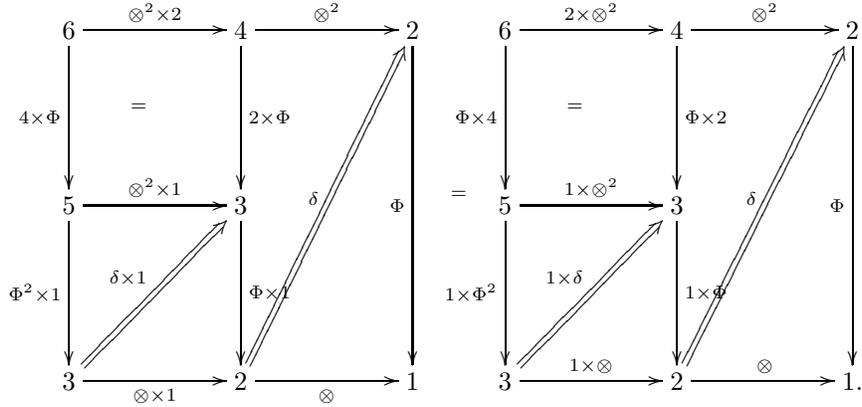
\xymatrix{  
6\ar[rr]^{\otimes^2\times 2}\ar[dd]_{4 \times \Phi}	&&4\ar[rr]^{\otimes^2}	\ar[dd]^{2 \times \Phi}		&&2  \ar[dddd]_\Phi
&6\ar[rr]^{2\times \otimes^2}\ar[dd]_{\Phi \times 4}	&&4\ar[rr]^{\otimes^2}	\ar[dd]^{\Phi \times 2}		&& 2 \ar[dddd]_\Phi
\\
\\
5\ar[rr]^{\otimes^2\times 1}\ar[dd]_{\Phi^2 \times  1}\ar@{}[uurr]^=&&3            \ar[dd]^{\Phi \times 1}		&&{}   \ar@{}[r]^=
&5\ar[rr]^{1\times \otimes^2}\ar[dd]_{1 \times \Phi^2}\ar@{}[uurr]^=	&&3             \ar[dd]^{1 \times \Phi }		&& {}            
\\
\\
3\ar[rr]_{\otimes \times 1} \ar@{=>}[rruu]^{\delta \times 1} 	&& 2\ar[rr]_{\otimes }\ar@{=>}[rruuuu]^{\delta}	                       		&&1                
&3 \ar[rr]^{1 \times \otimes  }\ar@{=>}[rruu]^{1 \times \delta}        &&2\ar[rr]^{\otimes  }\ar@{=>}[rruuuu]^{\delta}	&& 1.                
}$$
On the semantic level this means
$$\xymatrix{
(A \Phi A') \otimes (B\Phi B') \otimes (C\Phi  C') 
\ar[rr]^{\delta \otimes Id_{C\Phi C'}}
\ar[dd]_{Id_{A \Phi A'}\otimes \delta } &&
 (A\otimes B)\Phi (A'\otimes B')\otimes (C \Phi C') 
\ar[dd]^{\delta_{A\otimes B, A' \otimes B' , C, C'}} 
\\
\\
(A\Phi  A')\otimes (B\otimes  C) \Phi  (B'\otimes C') 
\ar[rr]^{\delta_{A,A',B \otimes C, B' \otimes C' }}  && 
 (A \otimes B \otimes C) \Phi  (A''\otimes  B'\otimes C')}$$
commutes. Setting the appropriate letters to $I$ (see page 
58 of \cite{Joy&Str}) gives us the braiding relations:  e.g.
$$\gamma_{A\otimes B, C} = (\gamma_{A,C} \otimes Id_B) \circ (Id_A \otimes \gamma_{B,C}).$$
\end{examp}

\begin{examp}
If we abandon the strictness (associativity) we get 
$$ \Tmon \KP \Tmon \cong \Tbraid . $$
The multiplications are $\otimes: 2 \longrightarrow 1$ and $\Phi :2 \longrightarrow 1$.
 Their respective reassociations are 
$\alpha : \otimes(\otimes \times 1) \Longrightarrow \otimes(1 \times \otimes)$ and 
 $\beta:\Phi(\Phi \times 1) \Longrightarrow \Phi(1 \times \Phi)$. 
Using similar results from last example, we can show that $A\otimes B \cong A\Phi B \cong
B\Phi A \cong B\otimes A$ as well as $\alpha \cong \beta$. The braiding relation is the only difference. 
Creatively pasting coherence conditions of $\delta$ we get several diagrams of the form:
\newline
$$\xymatrix{
9 \ar@{=}[r] \ar[dddd]_{(\Phi(\Phi \times 1))^3}^{\beta^3\Rightarrow}
&9\ar[rr]^{(1 \times \otimes)^3}\ar[dd]_{(1 \times \Phi)^3}	
&&
6\ar[rr]^{\otimes^3}	\ar[dd]_{(1 \times \Phi)^2}		&&
3  \ar[dd]_{1 \times \Phi} 
\\
\\
&6
\ar[rr]^{(1\times \otimes)^2 }\ar[dd]_{\Phi^3}
\ar@{=>}[uurr]^{1 \times \delta(\Phi,\otimes)}
&&
4
\ar[dd]_{\Phi^2} 
\ar[rr]^{\otimes^2} 
\ar@{=>}[rruu]^{\delta(1\times \Phi, \otimes)}
&&
2 \ar[dd]_{\Phi}       
\\
\\
3 \ar@{=}[r]
&3\ar[rr]_{1\times \otimes} 
\ar@{=>}[rruu]^{\delta(\Phi,1\times \otimes)}
 	&&
2\ar[rr]_{\otimes }
\ar@{=>}[rruu]^{\delta(\Phi,\otimes)}
&&
1                
\\
&3 \ar[rrrr]_{\otimes(\otimes \times 1)}^{\alpha \Uparrow}\ar@{=}[u]&&&&1\ar@{=}[u] .                
}$$
This diagram should be set equal to:
$$
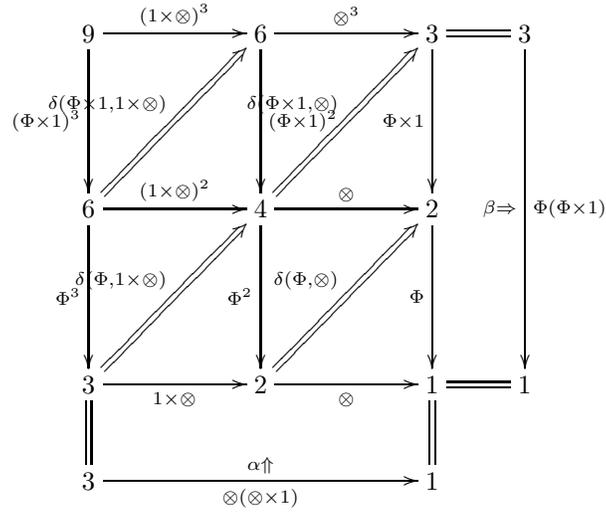
\xymatrix{
&9 \ar[rr]^{(1\times \otimes)^3}\ar[dd]_{(\Phi \times 1)^3}	&&
6\ar[rr]^{\otimes^3}	\ar[dd]^{(\Phi \times 1)^2}		&& 
3 \ar[dd]_{\Phi \times 1}\ar@{=}[r]&3 \ar[dddd]^{\Phi(\Phi \times 1)}_{\beta\Rightarrow}    
\\
\\
&6
\ar[rr]^{(1\times \otimes)^2}
\ar[dd]_{\Phi^3}
\ar@{=>}[uurr]^{\delta(\Phi \times 1,1 \times \otimes)}	
&&
4 
 \ar[dd]_{ \Phi^2 }\ar[rr]^{\otimes}
\ar@{=>}[rruu]^{\delta(\Phi \times 1,\otimes) }
 && 
2
 \ar[dd]_{\Phi}            
\\
\\
&3\ar[rr]_{1\times \otimes} 
\ar@{=>}[rruu]^{\delta(\Phi,1\times \otimes)} 
&&
2\ar[rr]_{\otimes }
\ar@{=>}[rruu]^{\delta(\Phi,\otimes)}	
  		&&
1\ar@{=}[r]&1
\\
&3 \ar@{=}[u]\ar[rrrr]_{\otimes(\otimes \times 1)}^{\alpha \Uparrow}&&&&1\ar@{=}[u]  .                
}$$

Setting 6 of the 9 variables to be $I$ (which will make some of the 2-cells into
the identity) and combining different forms of these diagrams will give us the famous dodecahedron:

$$
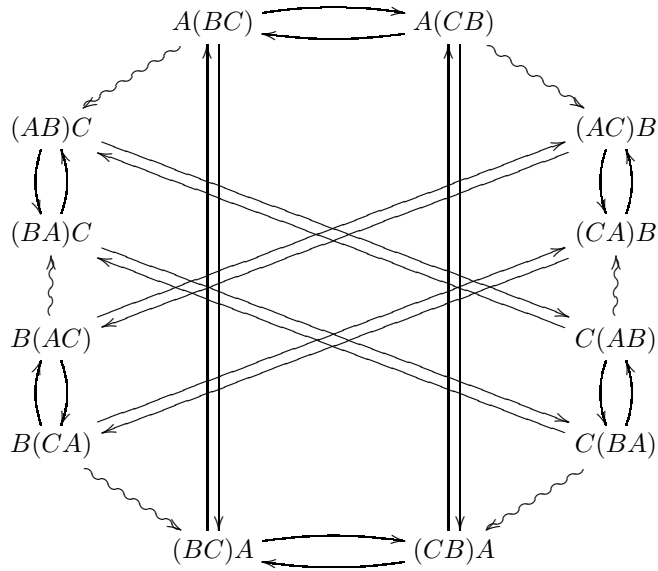
\xymatrix{
  &A(BC)\ar@/^/@{->}[rr] \ar@<.51ex>[ddddd] \ar@{~>}[dl] &     &A(CB)\ar@{~>}[dr]
\ar@/^/@{->}[ll] \ar@<.51ex>[ddddd]
\\
(AB)C\ar@/_/@{->}[d] \ar@<.5ex>[rrrrdd]&  &  &
&(AC)B\ar@/_/@{->}[d] \ar@<.5ex>[lllldd]
\\
(BA)C \ar@/_/@{->}[u]\ar@<.5ex>[rrrrdd]&   &  &  &(CA)B
\ar@/_/@{->}[u] \ar@<.5ex>[lllldd]
\\
B(AC)\ar@{~>}[u] \ar@/^/@{->}[d] \ar@<.5ex>[rrrruu]&   &  &  &C(AB)\ar@{~>}[u]
\ar@/_/@{->}[d] \ar@<.5ex>[lllluu]
\\
B(CA)\ar@{~>}[dr] \ar@/^/@{->}[u] \ar@<.5ex>[rrrruu]&  &  &
&C(BA)\ar@{~>}[dl]
\ar@/_/@{->}[u] \ar@<.5ex>[lllluu]
\\
&(BC)A\ar@/^/@{->}[rr] \ar@<.51ex>[uuuuu] &
&(CB)A\ar@/^/@{->}[ll] \ar@<.51ex>[uuuuu]
} $$

\vspace{2cm}

In order to distinguish the 
associativity isomorphisms from the commutativity isomorphisms
we draw associativity as 
$ \xymatrix{ \bullet \ar@{~>}[rr] & & \bullet }.$
The commuting dodecahedron is worth a few minutes of meditation.
The diagram actually indicates many equations.
Notice also that the rectangles commute  from that naturality of $\gamma$.
This naturality is, however, a semantical notion! $\gamma$ is merely a 2-cell 
in $\Tbraid$. The naturality comes from taking algebras in $\Cat$ where
2-cells are {\em natural transformations}. This leads us to ask what would 
happen if we took algebras in other 2-categories? Would the 
dodecahedron 
commute?

We would like to stress that the dodecahedron is not a ``new'' coherence structure. 
It is rather, a 2-dimensional statement that one reassociativity is a 
``homomorphism'' of the other reassociativity.
\end{examp}

\begin{examp}
Let $\Tsmon$ be the theory of strict monoidal categories with 
multiplication $\Phi:2 \longrightarrow 1$. Let $\Tsbraid$ be the
theory of strict braided monoidal categories with multiplication 
 $\otimes: 2 \longrightarrow 1$ and  braiding $\gamma: \otimes 
\Longrightarrow \otimes \circ tw$. Let $\T 0$ be the theory of pointed 
categories, that is, the theory of categories with a distinguished element to be thought of as 
a unit element. Then we have
$$ \Tsmon \KP \Tsbraid  \cong  \Tsbraid \KP \Tsmon \cong  \Tsym .$$
where $\Tsym$ is the theory of symmetric monoidal categories.

From the last two example, we know that $\Phi  \cong \otimes$ and that 
 $\delta$ can be made into a braiding. This braiding will be isomorphic to  $\gamma$.

The condition of the construction of $\KP$ shows us that 
$$\xymatrix{ 
	&2+2\ar[r]^{\Phi + \Phi}&1+1 \ar[dd]^{\otimes  tw}& =  & 	&2+2\ar[r]^{\Phi + \Phi}&
1+1\ddtwocell^{\otimes tw}_{\otimes}{^\gamma}  
\\
2+2 \ar[ru]^{\sigma}\dtwocell^{\otimes tw^2}_{\otimes^2}{^\gamma +\gamma} 
&&&& 
2+2 \ar[ru]^{\sigma}\ar[d]_{\otimes^2} & &
\\   
1+1 \ar[rr]_{\Phi} &{ }\ar@{=>}[ru]^{\delta(\otimes tw, \Phi)} & 1 & & 1+1 \ar[rr]_{\Phi} &{ }
\ar@{=>}[ru]^{\delta(\otimes , \Phi)} & 1}$$ 
Which translates into
$$\xymatrix{ 
A\otimes A' \otimes B \otimes B'
 \ar[rr]^{\gamma_{A,A'} \otimes \gamma_{B,B'} }
\ar[dd]_{Id \otimes \gamma_{A',B} \otimes Id}
 && A' \otimes A \otimes B' \otimes B 
\ar[dd]^{Id \otimes \gamma_{A,B'} \otimes Id}
\\
\\
A \otimes B \otimes A' \otimes B'
 \ar[rr]^{\gamma_{A\otimes B, A'\otimes B'}}
& & 
A' \otimes B' \otimes A \otimes B.}$$ commutes.

Setting $A=B'=I$ the unit (of both multiplications)
makes the top horizontal map and the right vertical map the identity. That leaves us with
$$ \gamma_{B,A'} \circ \gamma_{A',B}= Id$$ 
i.e.  symmetry.
\end{examp}

\begin{examp}
This is actually an example of something that does not work. 
Let $\Tsbraid$ be the theory of strict (associativity) braided monoidal
categories with multiplication $\otimes : 2 \longrightarrow 1$ and braiding
$\gamma: \otimes \Longrightarrow \otimes \circ tw$. Let $\TT_{Twist}$ be the 
theory with only one nontrivial generating 2-cell $\theta: Id_1 \Longrightarrow 
Id_1$ to be thought of as a twist of a ribbon. Let $\T 0 = \finbar$. One would expect 
that $\Tsbraid \KP \TT_{Twist}$ should be the theory of balanced categories
(see pg 65 of \cite{Joy&Str} or page 349 of \cite{Kassel} where they are  called ribbon
 categories. 
Truth be told, they assume a duality structure for the definition but it is not 
necessary for our needs.) The Kronecker product of these two theories forces 
the following equation

$$
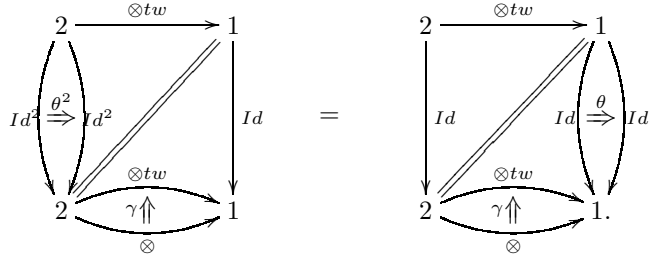
\xymatrix{
2\ar[rr]^{\otimes tw}\ddtwocell^{Id^2}_{Id^2}{^\theta^2}
&& 1\ar[dd]^{Id} && 2\ar[rr]^{\otimes tw}\ar[dd]^{Id} && 1\ddtwocell^{Id}_{Id}{^\theta}
\\
& &&=&&
\\
2 \rrtwocell_{\otimes}^{\otimes tw}{^\gamma} \ar@{=}[rruu]&&1&&2 \rrtwocell_{\otimes}^{\otimes tw}{^\gamma} \ar@{=}[rruu]&&1.
}$$

On the semantic level, this becomes 
\be  \gamma_{A,B}  \circ (\theta_A \otimes \theta_B) \quad =\quad  (\theta_{A\otimes B}) \circ \gamma_{A,B} \ee

This is very similar to the equation that is needed for a balanced category:
\be \gamma_{A,B}  \circ (\theta_A \otimes \theta_B)\circ \gamma_{B,A} \quad = \quad   \theta_{A\otimes B}. \ee
However these two equations are not the same! They would be the same if and only if the braiding were 
symmetric. Balanced tensor categories are ``part-way between braiding and
symmetry'' and it seems that the Kronecker product is too strong because it  makes the 
braiding symmetric. (See \cite{Paper2} for other structures  that are between braiding
and symmetry.) This (non)example
is strikingly similar to the 1-dimensional case where the Kronecker product of
the theory of monoids with the theory containing  one endomorphism of 1 (to be thought
of as the inverse) contains the theory of {\em commutative} (symmetric)  monoids.
\end{examp}

To what extent does the Kronecker product preserve the left quasi-adjoint
to $G^*$?  Consider $G_1:\T 1 \longrightarrow  \T 1'$ and  $G_2:\T 2 \longrightarrow  \T 2'$.
They induce the following diagram

$$\xymatrix{
\TA(\T 1 \KP \T 2, \Cat\ar@<1ex>[rr]^{Lan_{G_1 \otimes G_2}}_{\bot}
\ar[dd]_{\cong}
&&
\TA(\T 1' \KP \T 2', \Cat) \ar@<1ex>[ll]^{(G_1 \otimes G_2)^*} 
\ar[dd]^{\cong}
\\
\\
\TA(\T 1, \TA(\T 2, \Cat)) \ar@<1ex>[rr]^{Lan_{G_1}}_{\bot}
\ar[dd]_{\cong}
&&
\TA(\T 1', \TA(\T 2', \Cat)) \ar@<1ex>[ll]^{(G_1)^*} 
\ar[dd]^{\cong}
\\
\\
\TA(\T 2, \Cat) \ar@<1ex>[rr]^{Lan_{G_2}}_{\bot}
&&
\TA(\T 2', \Cat). \ar@<1ex>[ll]^{(G_2)^*} 
}$$

From the uniqueness of the quasi-adjoint of $(G_1 \otimes G_2)^*$
we may write
$$Lan_{G_1 \otimes G_2} \cong Lan_{G_1} \otimes Lan_{G_2}$$

Using this diagram, we can write new coherence results about Kronecker product theories from old coherence results.

Operads and theories are intimately related. They are two ways of 
describing algebraic structures on an object in a category.
Certain types of operads are in one-to-one correspondence
with 2-theories (see \cite{Thesis, Paper1} for a worked out example and \cite{Paper2} for a
general theory.) Markl \cite{Markl} has worked on a construct called a topological relative operad. It
is conjectured that this notion is nothing more then the operadic version of the Kronecker product. 

We would like to finish this paper by putting some of the facts that we have 
worked with in one commutative diagram. This diagram takes place in the
4 category of  ${\widetilde{\widetilde{\bf 3Cat}}}$. We shorten  the triple adjunction $c \vdash U \vdash
d \vdash \pi_0$ to 
$$
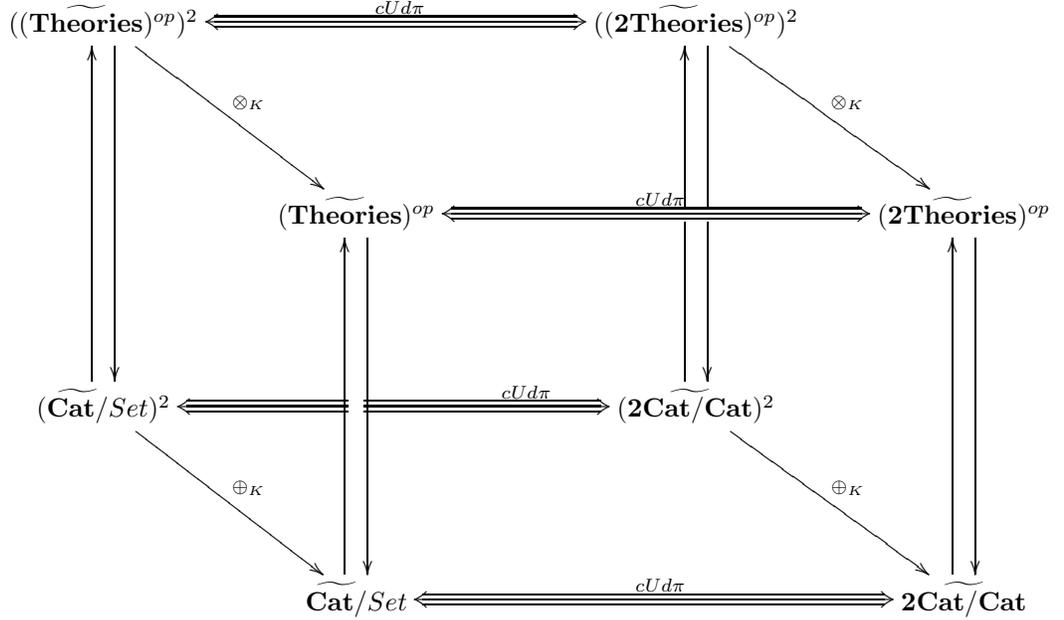
\xymatrix{
  \bullet \ar@3{<->}[rrr]^{cUd\pi}
&&& \bullet 
}$$
$$
\xymatrix{
((\Th )^{op})^2 
\ar[ddr]^{\otimes_K}
\ar@3{<->}[rrr]^{cUd\pi}
\ar@<1ex>[dddd]
&&& 
((\TTh)^{op})^2
\ar[ddr]^{\otimes_K}
\ar@<1ex>[dddd]|\hole
\\
\\
&(\Th)^{op}
\ar@3{<->}[rrr]^{cUd\pi}
\ar@<1ex>[dddd]
&&& 
(\TTh)^{op}
\ar@<1ex>[dddd]
\\
\\
(\Catset)^2 
\ar@3{<->}'[r][rrr]^{cUd\pi}
\ar[ddr]^{\oplus_K}
\ar@<1ex>[uuuu]
&&& 
(\TCatcat)^2
\ar[ddr]^{\oplus_K}
\ar@<1ex>[uuuu]|\hole |\hole
\\
\\
&\Catset \ar@3{<->}[rrr]^{cUd\pi}
\ar@<1ex>[uuuu]
&&& 
\TCatcat \ar@<1ex>[uuuu]
}$$

\begin{itemize}    
\item Top is syntax.
\item Bottom is semantics.
\item Left is one-dimensional universal algebra.
\item Right is two-dimensional universal algebra.
\item All diagonal maps are Kronecker products.
\end{itemize}
The fact that each of the the squares  commute was  either done in the paper or is left for the reader.

\section{Future directions}\

There  are many different directions in  which this work can be extended. 
An obvious generalization is multi-sorted 2-theories. More to the point, however,
would be 2-theories whose 0-cells are the free monoid on two generators $\lambda$
and $\rho$ corresponding to covariance  and contravarience. We may call 
such 2-theories ``bi-sorted 2-theories''.  Models/algebras of
such theories will be in a 2-category $\C$  that has both a product structure
and an involution $( ? )^{op}$. The prototypical example of such a category is  $\Cat$.
Algebras of these
theories would be functors that take $\lambda$ to $c$ and $\rho$ to $(c)^{op}$.
Using such a formalism, would help us understand the many structures  that 
demand contravarient functors. The list of structures that we could represent
 with such theories abound: monoidal closed categories, ribbon categories, traced
 monoidal categories, spherical categories etc. Algebraic functors and their left 
adjoints connecting all these structures would enlighten us about the relationship
 between them.

A further generalization of this paper would be monoidal 
2-theories. One can think of the our 2-theories as Cartesian 2-theories.
A monoidal 2-theory is similar to a Cartesian 2-theory but with a monoidal
product rather then a Cartesian product. Algebras will be (strict?) monoidal
preserving functors. This generalization would be of use to those who study 
k-linear categories with extra structure, relative coherence theory (see 
\cite{Paper5}) and quantum field theory (see next paragraph).

With the above two generalizations of this paper (and a healthy love of 
science fiction) we can apply bi-sorted monoidal 2-theories to the 
study of quantum field theory. Following Graeme Segal's conception
of conformal field theory, mathematical physicists have (see e.g.
\cite{Tillmann}) defined categories that look remarkably like 2-sorted
monoidal 2-theories. The 0-cells are finite families of circles oriented in
one of two ways ($\rho$ or $\lambda$). The 1-cells are to be thought of as 
``space-time segments'' from families of open circles to families of open circles.
The 2-cells are 
isotopy classes of diffeomorphisms that fix the boundary.
The 1-cells and 2-cells can have different structures depending on
what type of physical structure is of interest. The ``space-time
segments'' can be topological cobordisms (topological quantum field theory),
or Riemann manifolds (conformal field theory), or symplectic manifolds (symplectic
field theory). The tensor product in all of these theories is the disjoint union. 
There are many different functors between these 2-theories. For example there are
forgetful functors $U:{\bf T_{cft}} \lra {\bf T_{tqft}}$ and  
$U':{\bf T_{sft}} \lra {\bf T_{tqft}}$. What type of coherence results fall 
out of such 2-theory-morphisms? What does the ``free'' tqft for a 
given cft look like? What can we say about the (quasi?) adjoint functors 
induced from the inclusion of the $d$-dimensional tqft into the 
$d+1$-dimensional tqft? Is Tannaka duality \cite{Yetter}
nothing more then the reconstruction of the monoidal 2-theory
from its category of algebras? Is quantum field theory merely advanced 
universal algebra?

There are interesting questions arising from representation theory
Besides the triple adjunction $c \vdash U \vdash
d \vdash \pi_0$ between  ${\bf Theories}$ and 
$\widetilde{{\bf 2Theories}}$ there is yet another  
relationship between these two levels of structure 
that is less clear and needs to be studied. 
For every suitable algebraic (1-)theory $T$ and every $A \in Alg(T, Set)$
there is the
category of modules (suitably defined) for $A$. One of the main ideas 
in quantum groups is that the structure of the algebra $A$ is reflected
in the structure of the category of modules of $A$. Hence there is a
functor from ${\bf Theories}$ to  
 $\widetilde{{\bf 2Theories}}$  that takes $T$ to the 
2-theory of the structure of its category of modules.  For example, if $A$ is 
an old-fashioned algebra, then the category of modules is simply a
category. If
we add a coassociative comultiplication to $A$, then the category of modules 
inherits a strict monoidal structure. If the algebra has an involution 
(R-matrix, Drinfeld weak comultiplicaton structure, etc) then the category of
modules will have duality (braiding, monoidal structure, etc). Can this
functor from ${\bf Theories}$ to  
$\widetilde{{\bf 2Theories}}$ be formalized? Is there 
some type of inverse of this functor? Do we really
gain anything by going from the set with structure to the category  
with extra structure? Or can every theorem about categories with extra 
structure be understood  on the set with structure level? 
These constructions and questions are the
syntactical 
aspects of Tannaka duality.

Much work has been done lately to find the ``right'' definition
of a weak n-category. Allow me to give a definition of a weak
n-tuple category. A double category is a category object
in $\Cat$. Weakening this gives us a weak double category.
Iterating the construction of a double category gives 
us n-tuple categories. We are left asking what is a 
weak n-tuple category. Let $\TT$ be the 2-sketch of 
the theory of weak categories thought of as a set with endomorphisms 
and a partial operation. The partial operations make it a finite-limit
sketch rather then a 2-theory. We must extend the work done in section
4 to construct the Kronecker product of two finite-limit 2-sketches.
${\bf 2Alg(T, Set)}$ is the category of weak categories.
${\bf 2Alg(T \otimes T, Set) \cong 2Alg(T, 2Alg(T, Set))}$
is the 2-category of weak double categories.
${\bf 2Alg(T^{\otimes n}, Set)}$ is the category of weak n-tuple categories.
Coherence results will be induced by finite limit 2-sketch morphisms
of the form ${\bf T^{\otimes m} \lra T^{\otimes n}}$.

This paper does not close the door on functorial semantics.
There are many other
aspects of functorial semantics that we have not touched. For example, can we
characterize when a 2-category is a category of algebras for some 2-theory?
Can we characterize 2-functors as algebraic functors? When does a 2-theory
morphism $G$ induce a {\em  right} (quasi-) adjoint to $G^*$? etc.

Further study needs to be done on the intimate relationship between
 2-theories and 2-monads (see e.g. Blackwell et al \cite{Blackwell} 
and Lawvere \cite{Lawvereb} ).
 The study of the
 connection between theories and monads spawned much insight into
 both structures and we are
 sure that the same study of their two-dimensional analogs would be 
just as fruitful.

Computer science has long since coopted algebraic theories for its own use. 
Wagner \cite{Wagner} is a survey article of all these types of theories (e.g. ordered
 theories, iteration theories, rational theories, iterative theories etc.) Such 
generalizations have been used in diverse fields of computer science such as 
context-free grammars, flowchart semantics, recursion schemata and 
recursively defined domains. There 
is surely room to do similar generalizations for 2-theories. There are 
other areas of computer science that would benefit 
from  a study of 2-theories. 
Seely \cite{Seely} has an approach to lambda-calculus and computation using 
2-categories. The entire area of linear logic uses categories with structure that 
could and should be put into a 2-theoretic context.

There is, obviously, a deep connection between higher dimensional category theory and  
homotopy theory. However, it is not too obvious what the connection actually is. 
Perhaps we would be able to better understand this connection by looking at the algebra
case. Over the past few years there has been a tremendous amount of work 
``homotopy algebras'' or ``deformation algebras'' ( or  `` $( -)_{\infty}$
algebras''). These might all be formulated using theories. They all 
have some connections to homotopy. We conjecture that with all these algebras, 
their category of modules have added structure that can be formalized with
a 2-theory. So, in a sense, the homotopy aspects of these algebras can be seen
in their 2-theoretic formulation. Can more be said on this topic?

Department of Computer and Information Science\\
Brooklyn College, CUNY\\
Brooklyn, N.Y. 11210\\
email: noson@sci.brooklyn.cuny.edu
\end{document}